\documentclass[preprint,12pt]{elsarticle}

\usepackage{amssymb}
\usepackage{amsmath}
\usepackage{enumitem}
\usepackage{graphicx}
\usepackage[hidelinks]{hyperref}
\usepackage{makecell}

\journal{Journal of Sound and Vibration}

\begin{document}
	
	\begin{frontmatter}
		
		
		
		\title{Dynamical integrity estimation in time delayed systems: a rapid iterative algorithm}
		
		
		\author[label1,label2]{Bence Szaksz}
		
		\author[label1,label3]{Gabor Stepan}
		
		\author[label1,label2]{Giuseppe Habib}
		
		\affiliation[label1]{organization={ Department of Applied Mechanics, Faculty of Mechanical Engineering, Budapest University of Technology and Economics},
			city={Budapest},
			postcode={1111}, 
			country={Hungary}}
		
		\affiliation[label2]{organization={MTA-BME Lendület “Momentum” Global Dynamics Research Group, Budapest University of Technology and Economics},
			city={Budapest},
			postcode={1111}, 
			country={Hungary}}
		
		\affiliation[label3]{organization={ELKH-BME Dynamics of Machines Research Group},
			city={Budapest},
			postcode={1111}, 
			country={Hungary}}
		
		\begin{abstract}
			The robustness of dynamical systems against external perturbations is crucial in engineering; however, it is often overlooked for the lack of methods for rapidly computing it. This paper proposes a novel algorithm for estimating the robustness of systems subject to time delay. More precisely, the algorithm iteratively estimates the so-called local integrity measure (LIM), that is, the radius of the largest hypersphere centered at the fixed point and located within its basin of attraction. Since time-delayed systems are infinite dimensional, initial conditions are restricted to a constrained type of initial functions. The semi-discretization method is used for rapidly simulating the dynamics of the systems, while trajectories are classified as converging or diverging using a subdivision of the reduced phase space into cells. The algorithm was tested on four different mechanical systems, and in all cases it very quickly provided an accurate estimation of the LIM. Moreover, it enabled the study of LIM trends in a multi-dimensional parameter space, which would have been unfeasible with alternative methods. This breakthrough in computational efficiency has important implications for engineering design, allowing for careful consideration of dynamical integrity and enhancing the safety and reliability of engineered systems, especially in the presence of time delays.
		\end{abstract}
		
		
		
		
	\end{frontmatter}
	
	
	
	\section{Introduction}
	
	In engineering practice, dynamical systems are often analyzed through linearization, which has several benefits. Linear systems have only one equilibrium, and their dynamics can be investigated with straightforward methods. Nowadays, as industries develop increasingly complex and sensitive products,  understanding a system's nonlinear behavior is becoming more and more important. For example, to reduce weight and energetic costs, modern structures are often slender, which amplifies nonlinear effects.
	
	Nonlinear dynamical systems, unlike linear ones, may have several equilibria and/or limit cycles, multiple of which can be stable simultaneously. Bistability is common in many systems and each stable equilibrium has its own basin of attraction (BoA). Large perturbations can cause the system to transition from the desired equilibrium to an undesired one, leading to failures and accidents, such as shimmy in vehicle wheels \cite{pacejka1965analysis,beregi2019bifurcation,horvath2022stability,habib2023towed}, flutter in airplane wings \cite{lind2003flight,verstraelen2017experimental,basta2021flutter,takarics2020active,drachinsky2022flutter}, turbulent flows \cite{kerswell2018nonlinear,cherubini2015nonlinear}, traffic jams \cite{nagatani2002physics,orosz2010traffic,molnar2021delayed}, 
	robot control \cite{veraszto2017nonlinear,szaksz2022transient,habib2022bistability,bartfai2022hopf}, machining operations \cite{dombovari2015bistable,molnar2019closed,iklodi2022bi}, electric blank-outs \cite{pourbeik2006anatomy,gajduk2014stability}, human balance \cite{smith2017basins,zakynthinaki2010modeling} and predator-prey ecosystems \cite{aguirre2014bifurcations,arancibia2019basins} to mention a few.
	
	Several methods exist in the literature for investigating the robustness of a desired equilibrium.  The most commonly implemented method consists of performing numerical simulations from a grid of initial conditions, which enables the detection of the basin of attraction. However, this method is intrinsically inefficient and computationally expensive. It is used only to obtain two-dimensional sections of the basin and it is prohibitive for high-dimensional systems \cite{rega2021global}. Probabilistic approaches mainly based on Monte Carlo simulations are more efficient \cite{martiniani2016structural,sprott2015classifying,yan2021statistical}, but their outcomes are not compatible with integrity measures \cite{lenci2019global}. Analytical methods are usually based on Lyapunov functions, with the help of which a subset of the basin of attraction can be determined \cite{ratschan2010providing,grinberg2015boundary,biemond2014estimation}. However, generating an appropriate Lyapunov function can be challenging and is nearly impossible for systems with many degrees-of-freedom (DoF).
	
	Cell mapping, introduced by Hsu, is a highly effective method for estimating the basin of attraction \cite{hsu1980theory,hsu1986cell,hsu2013cell}. This method discretizes the phase space into cells and calculates short trajectories from each cell to determine the flow's tendency. Cell mapping is suitable for parallel computation and is efficient for relatively large systems with many degrees of freedom. Seeing its potential, the researchers developed the algorithm further in the past years \cite{sun2018cell,liu2016global,andonovski2020six}. However, this method requires to store information corresponding to each cell, which yields memory issues in case of large systems.
	
	A different approach is proposed in \cite{habib2021dynamical}, which directly looks for the dynamical integrity without computing the BoA.  Specifically, it estimates the local integrity measure (LIM) of nonlinear dynamical systems, which is the radius of the largest hypersphere entirely within the BoA and centered at the desired stable fixed point. This paper aims to generalize this algorithm to estimate the LIM in systems subject to time delay.
	
	Time delay is present in every controlled dynamical system, regardless of whether it is human or robot controlled since the data processing and actuation requires time \cite{hu2003dynamics,szaksz2022delay,szaksz2022nonlinear}; moreover, delay also arises when a system has some kind of memory effect \cite{sun2017effect,de2020memory,stepan2001modelling}. Delay differential equations (DDEs) describe the dynamics of systems subject to time delay, such as the above mentioned traffic jams, wheel shimmy, and machine tool vibrations \cite{beregi2019bifurcation,orosz2010traffic,szaksz2023delay,altintas2020chatter,stepan2014cylindrical,kiss2022control}. 
	Investigating DDEs is challenging because they have an infinite-dimensional phase space requiring a function of time as initial condition \cite{hale2013introduction,stepan1989retarded}. Therefore, the definition of the basin of attraction of a fixed point should be reconsidered as well.
	
	Because of the above mentioned difficulties, the BoA of time delayed systems are very rarely computed. Analytical methods are usually based on Lyapunov-Krasovskii functionals but their calculation is complicated even for medium size systems \cite{biemond2014estimation}. Another option is to restrict the initial condition to a specific type, as constant, linear, jump or free vibration initial functions \cite{yoshida2021basins,scholl2020norm}.  Furthermore, there are also studies in the literature about how to obtain BoA from human balance experiments, when the relevant human reaction time is also considered \cite{zakynthinaki2010modeling,smith2017basins}.
	
	This paper suggests a novel algorithm for the fast and accurate estimation of the LIM in time delayed systems. Semi-discretization is applied for the time-integration of trajectories corresponding to different initial conditions; the estimated value of the LIM is iteratively reduced if a diverging trajectory is found within the currently assumed
	hypersphere of convergence. To speed up the algorithm, it is 
	continuously monitored whether a new trajectory is converging to a previously categorized one.
	
	The paper is organized as follows: Sec. \ref{sec:definitions} defines the basin of attraction for time-delayed systems and introduces the corresponding LIM and a distance definition to compare states with different units. Sec. \ref{sec:algorithm} describes the proposed algorithm, followed by four case studies in Sec. \ref{sec:casestudies}, which show the algorithm's potential to provide fast and accurate estimations of the systems' dynamical integrity, which is useful from an engineering standpoint. Finally, Sec. \ref{sec:discussion} discusses the advantages and limitations of the proposed algorithm, and provides concluding remarks.
	
	\section{Definitions}\label{sec:definitions}
	\subsection{Basin of attraction}
	
	The BoA is defined as the set of initial conditions that converge to a particular equilibrium point or to a limit cycle in a dynamical system. In the case of delay differential equations (DDEs), the phase space is infinite-dimensional \cite{michiels2007stability}, which makes it difficult to identify the BoA for arbitrary initial functions. Therefore, the BoA is usually defined for specific sublevelsets in the space of initial functions \cite{biemond2014estimation,leng2016basin,yan2021statistical}, such as constant or linear initial functions, or for initial conditions where the system performs free vibrations \cite{yoshida2021basins,scholl2020norm}; in the current study, primarily the latter one is applied, as discussed later. 
	
	Here, we define the BoA based on the headpoint of the constrained initial condition. A point in phase space belongs to the BoA if it is the headpoint of an initial condition that yields a converging solution. Note that, since the BoA is constructed based on the initial condition headpoints (ICHs), there may be diverging trajectories the ICH of which is outside of the BoA, while the trajectory itself enters and then leaves the BoA permanently. 
	
	\begin{figure}
		\centering
		\includegraphics{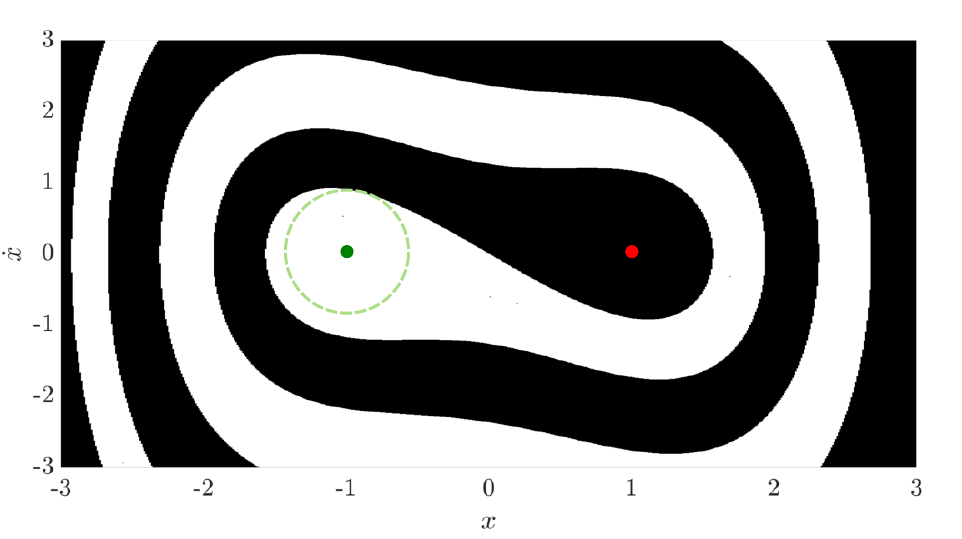}
		\caption{Basins of attraction of the two stable fixed points of the delayed Duffing equation $\ddot{x}(t)+0.2\dot{x}(t)-x(t-0.1)+x^3(t) = 0$. Green and red dots refer to the desired and to the undesired fixed points, respectively, while the dashed green curve is the hypersphere of convergence the radius of which is the LIM.}
		\label{fig:DuffingBoA}
	\end{figure}
	
	\subsection{Local integrity measure}
	
	The BoA can be a intricate subspace, sometimes with intermingled or fractal-like boundaries \cite{soliman1989integrity}.  The computation of the BoA is very expensive, while from an engineering point of view, the intermingled and fractal areas should be avoided since in those domains small perturbations are enough to leave the BoA. Dynamical integrity measures were introduced to overcome these issues.
	
	Different types of dynamical integrity measures have been proposed in the literature, such as, the global integrity measure, the impulsive integrity measure,  the stochastic integrity measure, the integrity factor or the actual global integrity measure \cite{thompson1989chaotic,lenci2003optimal,rega2005identifying}. Although they are diverse, they all aim at quantifying the dynamical integrity of the system, and they exhibit similar trends for varying system parameters \cite{lenci2019global}.
	To assess the system's robustness against perturbations in initial conditions, we use the LIM  \cite{soliman1989integrity}, which is defined as the radius of the largest hypersphere centered at the fixed point that lies entirely within the BoA as it was mentioned in the introduction. It is a conservative, simple measure, that is easy to estimate in an iterative manner.
	
	Figure~\ref{fig:DuffingBoA} shows the basins of attraction of the two stable equilibria of the delayed Duffing equation: $\ddot{x}(t)+0.2\dot{x}(t)-x(t-0.1)+x^3(t) = 0$. Considering the equilibrium point $(-1,0)$ as the desired one, the green dashed circle denotes the hypersphere of convergence, the radius of which is the LIM. 
	
	\subsection{Distance metric}
	Usually, the coordinates of the phase space have diverse units, therefore, there is a need for a distance definition, where the different directions are taken into account with predefined weights. These weights can be chosen based on the degree of sensitivity of the system in different directions. For example, in an engineering problem, if the system is more sensitive to perturbations in velocity rather than position, then larger weights can be assigned to velocity coordinates.
	
	Following the suggession of \cite{habib2021dynamical}, we define the distance metric based on the mechanical energy of the undamped linear system without control. This approach is advantageous as it takes into account the dynamics of the system and provides a physically relevant measure.
	
	Let us consider the $n$ DoF dynamical system in the form:
	\begin{equation}   \mathbf{M}\ddot{\mathbf{x}}+\mathbf{C}\dot{\mathbf{x}}+\mathbf{K}\mathbf{x}=\mathbf{0}\,,
	\end{equation}
	where $\mathbf{M}$, $\mathbf{C}$ and $\mathbf{K}$ are the mass, damping and stiffness matrices, respectively, while $\mathbf{x}$ is the state vector of the physical coordinates. 
	Neglecting the damping, the system can be rewritten in the space of the modal coordinates $q_i$ ($i = 1,2,\dots,n$) as
	\begin{equation}
	\ddot{q}_i+\omega^2_i q_i=0,\quad \mathrm{for}\quad i = 1,2,\dots,n\,,
	\end{equation}
	where $\omega_i$ is the natural angular frequency of the $i$-th mode and $n$ is the degree of freedom. The normalized mechanical energy stored in the $i$-th mode is $\omega_i^2 q_i^2+\dot{q}_i^2$; based on which, let us introduce the definition of the distance as:
	\begin{equation}\label{eq:distancedef}
	d=\sqrt{\sum_{i=1}^n(\omega_i^2 q_i^2+\dot{q}_i^2)}\,.
	\end{equation}
	Thus, the distance between two points $\mathrm{A}$ and $\mathrm{B}$ of the phase space takes the form
	\begin{equation}\label{eq:distance}
	d_\mathrm{AB}=\sqrt{\sum_{i=1}^n \left(\alpha_i(q_{i\mathrm{A}}-q_{i\mathrm{B}})^2\right) + \sum_{i=1}^n \left(\alpha_{n+i}(\dot{q}_{i\mathrm{A}}-\dot{q}_{i\mathrm{B}})^2\right)}
	\end{equation}
	with the weight vector $\alpha$ that assumes the form
	\begin{equation}\label{eq:distanceweights}
	\alpha = [\omega_1^2, \omega_2^2, \dots, \omega_n^2,1,\dots,1]\,.
	\end{equation}
	
	Of course, in many cases the modal decomposition cannot be carried out, such as for non-mechanical systems or for odd-dimensional systems for example. In that case the simplest way is to set all the weights to one or tune them based on some engineering intuition. We note that, in general, the choice of the distance weights has a minimal effect on the trend of the dynamical integrity in the parameter space. Besides, the working principle of the algorithm for the LIM estimation is not affected by the weight values.
	
	\section{Proposed algorithm}\label{sec:algorithm}
	
	\subsection{Overview of the algorithm}\label{sec:iterativeproc}
	This paper suggests a novel algorithm for the fast and accurate estimation of the LIM in time delayed systems. After defining the space boundary and some other parameters, the iteration can be started. The algorithm is based on a semi-discretization scheme \cite{insperger2011semi}, which calculates the trajectories corresponding to a constrained initial condition. Each trajectory is categorised whether it is converging to the desired fixed point or not. To speed up the algorithm, the method continuously checks whether a new trajectory is converging to a previously categorized one. If a divergent trajectory is found, the estimated value of the LIM is reduced accordingly, and the area of interest is reduced to the corresponding hypersphere of convergence. Three options are suggested for the selection of the headpoint of the next initial condition: random selection, bisection method, and that point of the last diverging trajectory, which is closest to the desired fixed point; these methods can also be combined. 
	
	\subsection{Modal decomposition}
	The first step of the algorithm is to perform a modal decomposition of the dynamical system; this step provides the natural angular frequencies that determine the actual weight vector. The algorithm then transforms the nonlinear dynamics into the space of the modal coordinates, which has the advantage that the free vibration initial conditions are easy to handle with mode shapes. If the dimension of the system is odd or the system does not possess vibration modes, this step is skipped and the algorithm works with user-defined weights in the space of the physical coordinates.
	
	\subsection{Semi-discretization}
	
	Time series of the system under study are obtained through the semi-discretization method \cite{insperger2011semi}.
	The semi-discretization method was originally introduced for linear DDEs, and the key idea is to consider the delayed terms to be piece-wise constant. In this way, the DDE can be rewritten as a non-homogeneous ordinary differential equation (ODE) over the sampling time intervals, which can be solved piece-wise with the help of Duhamel's (variation of coefficients) formula. The approach can be implemented also to nonlinear systems \cite{insperger2011semi,molnar2018bifurcation}, as explained below.
	
	Let us consider the first order nonlinear DDE
	\begin{equation}\label{eq:semibasiceq}
	\dot{\mathbf{y}}(t)=\mathbf{A}\mathbf{y}(t)+\mathbf{B}\mathbf{y}(t-\tau)+\mathbf{g}(\mathbf{y}(t),\mathbf{y}(t-\tau))\,,
	\end{equation}
	where $\mathbf{y}$ is the state vector, $\mathbf{A}$ and $\mathbf{B}$ are the coefficient matrices of the non-delayed, and of the delayed states, respectively, while the function $\mathbf{g}$ contains the nonlinear terms. 
	Let us assume a discretization with time step $h\ll\tau$, introduce the sampling instants $t_i=ih$, and approximate both the delayed terms and the nonlinear terms with a piece-wise constant function over the time interval $t\in[t_i,t_{i+1})$ leading to
	\begin{equation}\label{eq:semipiecewise}
	\dot{\mathbf{y}}(t)=\mathbf{A}\mathbf{y}(t)+\mathbf{B}\mathbf{y}(t_{i-r})+\mathbf{g}(\mathbf{y}(t_i),\mathbf{y}(t_{i-r}))\,,
	\end{equation}
	where  $r$ is a positive integer called sampling delay number.
	
	The effect of the delayed terms appears first with a delay of $r h$ and does not change until the next sampling; this yields a saw-tooth like time-varying delay:
	\begin{equation}
	\rho(t)=rh + t - t_i\,,\quad t\in[t_i,t_{i+1})\,.
	\end{equation}
	The corresponding average delay is 
	\begin{equation}\label{eq:tauaver}
	\bar{\tau}=\left(r+\frac{1}{2}\right)h\,.
	\end{equation}
	
	\begin{figure}
		\centering
		\includegraphics[width=\linewidth]{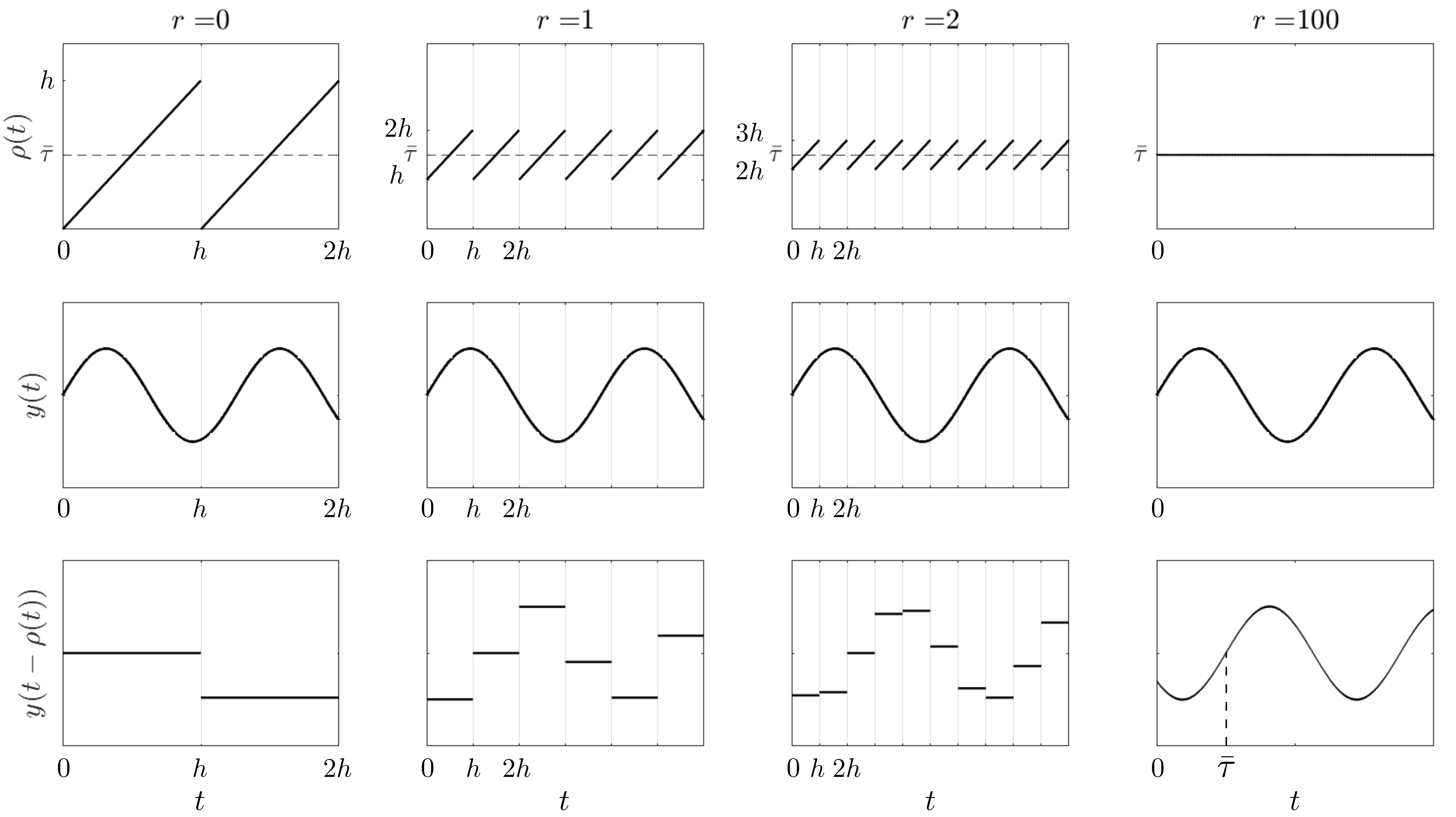}
		\caption{Effect of sampling of the sinusoidal signal $y(t)$ for various values of the sampling delay number $r$. The panels in the first row present the time-varying time delay $\rho(t)$, the panels in the middle show the original signal, while the panels in last row refer to the sampled signal with average time delay $\bar{\tau}$.}
		\label{fig:semidisc}
	\end{figure}
	
	If $r\rightarrow \infty$ and $h \rightarrow 0$ such that Eq.~\eqref{eq:tauaver} holds, the solution of the discretized DDE tends to the solution of (6) with $\tau=\bar{\tau}$ as it is visualized in Fig.~\ref{fig:semidisc}.
	
	Solving \eqref{eq:semipiecewise} with an initial condition $y_i$ provides the initial condition for the next $t\in[t_{i+1},t_{i+2})$ interval as
	\begin{equation}
	\mathbf{y}_{i+1}=\mathbf{e}^{\mathbf{A}h}\mathbf{y}_i+ \int\limits_0^h \mathbf{e}^{\mathbf{A}(h-s)}\mathrm{d}s (\mathbf{B} \mathbf{y}_{i-r}+\mathbf{g}(\mathbf{y}_i,\mathbf{y}_{i-r}))\,.
	\end{equation}
	
	Thus, the dynamics of the original nonlinear delayed system \eqref{eq:semibasiceq} can be approximated with the discrete map
	\begin{equation} \label{eq:semidiscmap}
	\begin{bmatrix}
	\mathbf{y}_{i+1}\\
	\mathbf{y}_{i}\\
	\mathbf{y}_{i-1}\\
	\vdots \\
	\mathbf{y}_{i-r+1}
	\end{bmatrix} = 
	\begin{bmatrix}
	\mathbf{P}&\mathbf{0}&\dots&\mathbf{0}&\mathbf{Q}\mathbf{B}\\
	\mathbf{I}&\mathbf{0}&\dots&\mathbf{0}&\mathbf{0}\\
	\mathbf{0}&\mathbf{I}&\dots&\mathbf{0}&\mathbf{0}\\
	\vdots&&\dots&&\vdots\\
	\mathbf{0}&\mathbf{0}&\dots&\mathbf{I}&\mathbf{0}
	\end{bmatrix}
	\begin{bmatrix}
	\mathbf{y}_{i}\\
	\mathbf{y}_{i-1}\\
	\mathbf{y}_{i-2}\\
	\vdots \\
	\mathbf{y}_{i-r}
	\end{bmatrix}
	+
	\begin{bmatrix}
	\mathbf{Q}\mathbf{g}(\mathbf{y}_i,\mathbf{y}_{i-r})\\
	\mathbf{0}\\
	\mathbf{0}\\
	\vdots \\
	\mathbf{0}
	\end{bmatrix}\,,
	\end{equation}
	where
	
	\begin{equation}\label{eq:semidiscmatr}
	\mathbf{P}=\mathbf{e}^{\mathbf{A}h},\quad \mathrm{and}\quad 
	\mathbf{Q}= \int\limits_0^h \mathbf{e}^{\mathbf{A}(h-s)}\mathrm{d}s  \,.
	\end{equation}
	Mathematical proof of the convergence of the discretized system to the original one is provided in \cite{insperger2011semi}. The main advantage of this formulation is that it is significantly faster than numerically integrating in time the original DDE to ensure the same level of accuracy.
	
	\subsection{Initial conditions}
	
	The above described semi-discretization based mapping requires an initial condition for the $t\in[-\tau,0]$ interval, which is sampled with the time step $h$. Since the convergence of the solution depends on the entire initial function, and not only on its headpoint, a constrained type of the initial function is considered to obtain a consistent BoA \cite{scholl2020norm,yoshida2021basins}. Each type of initial function corresponds to a specific physical situation. For example, a dynamical system may lose stability because of an impulse-like perturbation, or if the delayed control is turned on during the free vibration of the system. In this section, four types of initial conditions are discussed, although, the algorithm can handle any kind of sampled initial function:
	
	\begin{figure}
		\centering    \includegraphics[width=0.9\linewidth]{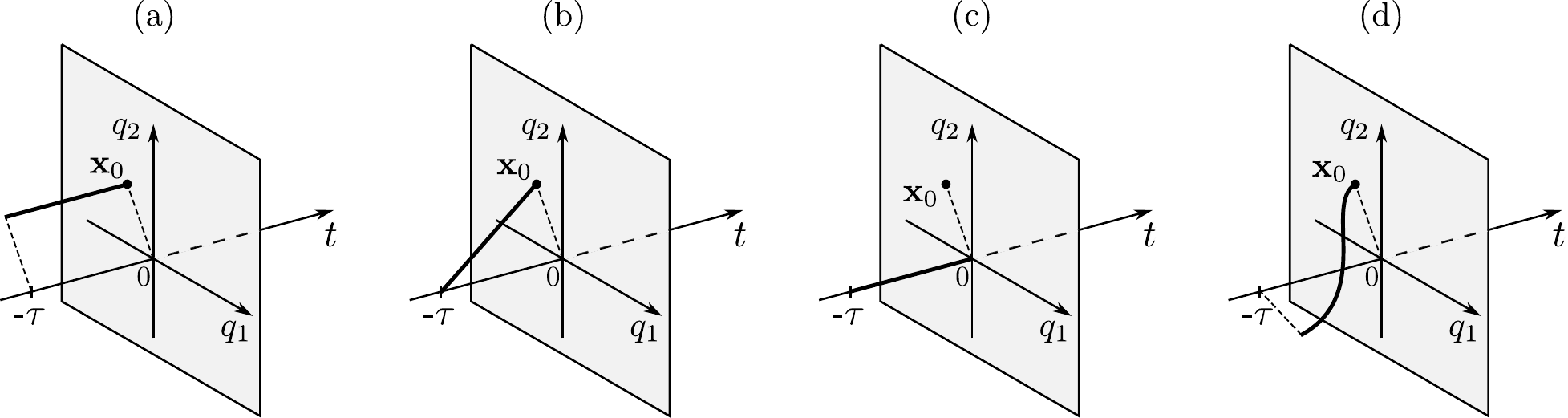}
		\caption{Different types of initial functions with the same headpoint $\mathbf{x}_0$: (a) constant, (b) linear, (c) jump, (d) free vibration initial function.}
		\label{fig:initconds}
	\end{figure}
	
	\begin{enumerate}[label=(\alph*)]
		\item The simplest type of initial function is to choose a \textit{constant} initial condition in each coordinate (see Fig.~\ref{fig:initconds}(a)), which was used in several studies \cite{scholl2020norm,yoshida2021basins}. Indeed, it is easy to implement and provides useful results, but it leads to a contradiction in real physical systems since the position and the velocity cannot be a nonzero constant at the same time. 
		
		\item Another option is to apply a \textit{linear} initial condition between the equilibrium point and a specified point $\mathbf{x}_0$ in the phase space \cite{scholl2020norm}; see Fig.~\ref{fig:initconds}(b). However, this leads again to a contradiction considering the physical relation between the positions and the velocities.
		
		\item In case of impulse-like perturbations, the system is in equilibrium state after which the velocity jumps to a certain value that is referred to by $\mathbf{x}_0$ in Fig.~\ref{fig:initconds}(c) \cite{scholl2020norm,yoshida2021basins}. This provides physically relevant solutions, since this \textit{jump} type of initial function satisfies the original DDE yielding special type of integro-differential equations \cite{burton2005volterra,burton2016existence}.
		
		\item During the case studies, we mainly use a \textit{free vibration} based initial condition, when the undamped linear harmonic oscillation is considered without control; see Fig.~\ref{fig:initconds}(d). First, a point $\mathbf{x}_0$ in the phase space is selected as the headpoint of the initial function, then the corresponding history is calculated based on the free vibration of the system. Although the nonlinearity and damping are neglected, the initial condition is still related to the dynamics of the system. We assume that the delayed control turns on after $t=0$.
	\end{enumerate}
	
	Despite the intrinsic differences of the discussed types of initial functions, we note that they do not qualitatively affect the results provided by the algorithm, as it will be illustrated later.
	
	\subsection{Trajectory characterization}
	Before starting the iteration, the operator has to define the location of the desired fixed point and the space boundary (a hyperrectangle), within which the LIM is looked for. Furthermore, cell subdivision is applied to the area of interest; for simplicity, the same number of cells is used for each coordinate direction. At the sampling instants, not only the exact point of the phase space is saved, but also the cell to which the point belongs. Therefore, the trajectory can be represented as a series of cells, which allows for efficient characterization as it is explained below.
	
	\begin{figure}
		\centering
		\includegraphics[width=\linewidth]{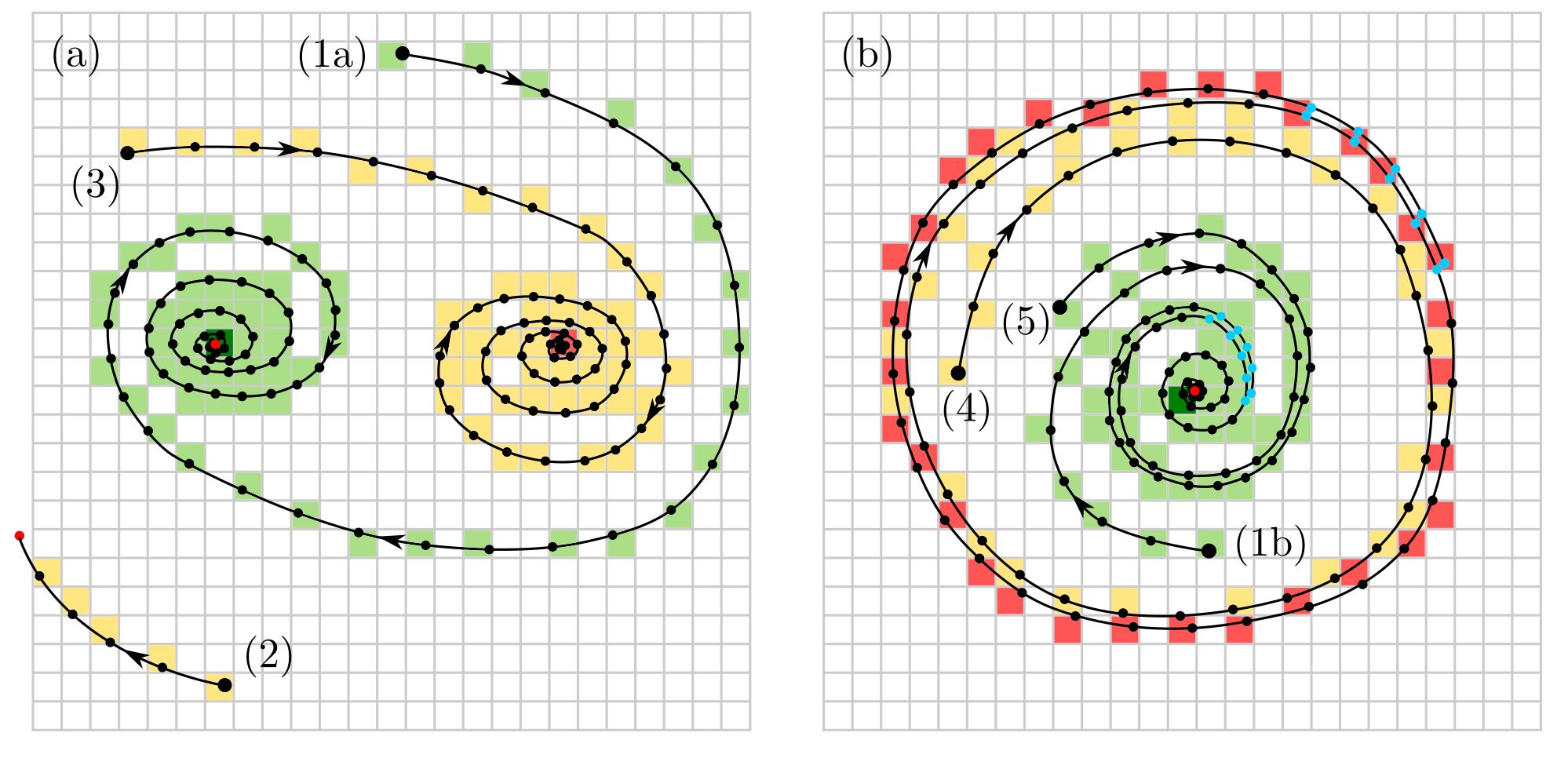}
		\caption{Trajectory characterization in two different types of dynamical systems. The cells corresponding to the sampled converging and diverging trajectories (with respect to the desired fixed point) are colored green and yellow,  respectively. Blue dots refer to a series of data points where the trajectory is found convergent to itself or to a previously categorized trajectory.}
		\label{fig:categorisation}
	\end{figure}
	
	Figure~\ref{fig:categorisation} shows typical trajectories in different types of dynamical systems. These trajectories can be classified into two main categories with some subcategories:
	\begin{enumerate}[label=(\alph*)]
		\item \textit{Converging to the desired solution} 
		\begin{itemize}    
			\item The desired equilibrium is defined in the pre-processing step. The mapping stops and the corresponding trajectory is characterized as convergent if the  solution remains within the cell of the equilibrium for a time larger than $\tau$ (see trajectories (1a) and (1b) of Fig.~\ref{fig:categorisation}).  Because of the exponential nature of the solutions, as they get closer to the fixed point, the settling slows down and they may oscillate for a long time around the equilibrium before entering the corresponding cell. Therefore, in cases of small damping, it is useful to define not only the cell of the equilibrium point as convergent, but also the neighbouring cells.
		\end{itemize}
		\item \textit{Non-converging to the desired fixed point}
		\begin{itemize}    
			\item \textit{Diverging out of the space boundary.} If the trajectory crosses the predefined space boundary, than it is categorized as divergent; see trajectory (2) of Fig.~\ref{fig:categorisation}(a).
			\item \textit{Converging to a new, previously unknown fixed point.} Nonlinear systems may have multiple coexisting equilibrium points. If the trajectory spends time more than $\tau$ within a cell, then that cell is classified as a cell containing a new fixed point, and the trajectory is categorized as divergent; see trajectory (3) of Fig.~\ref{fig:categorisation}(a). If it occurs that the trajectory crosses a cell very slowly, dwelling there for a time longer than $\tau$ before leaving it, a non-existing (`ghost') fixed point might erroneously be identified. In order to avoid this issue, the trajectory may be required to spend time more than $\bar{t}>\tau$ within the same cell after which the cell is considered to contain a new fixed point.
			\item \textit{Converging to a periodic solution.} If the trajectory converges to a periodic solution, it is still categorized as divergent. Since the latest subsequent cells with length $\tau$ describe the actual state of the solution, the algorithm continuously monitors whether this sequence of cells are found in the earlier part of the discrete trajectory. If so, then the time-integration is interrupted and the trajectory is categorized as a periodic solution. For example, in Fig.~\ref{fig:categorisation}(b), trajectory (4) converges to a periodic orbit, which is detected when a sequence of five cells is repeated within the same trajectory (denoted by blue dots). 
			
			In the case of small damping, trajectories slowly spiraling towards an equilibrium point might be classified erroneously as periodic solutions (a problem typically encountered with the cell mapping method as well \cite{vio2007bifurcation}). To avoid this event, the number of times a series of cells must be repeated in a trajectory to classify it as periodic should be increased.
		\end{itemize}
	\end{enumerate}
	
	Suppose that a given number of subsequent trajectory points passes through the same cells of a previously classified trajectory. In that case, the current simulation is interrupted, and the new trajectory is classified as the one reached (see trajectory (5) of Fig.~\ref{fig:categorisation}(b)). Indeed, this event is rather common, as trajectories tend to approach an invariant manifold after an initial transient \cite{cirillo2016spectral,cenedese2022data}. Accordingly, this strategy enables the significant reduction of computational time.
	
	Furthermore, there are some cases when the solution cannot be categorized for a long time. To avoid too long simulations, the semi-discretization map is stopped when the time reaches a predefined value. In this case, either a human operator should categorize the trajectory based on a phase space representation, or the trajectory is classified as diverging, which has the engineering relevance that both the diverging and the slowly converging solutions are undesired. 
	
	To obtain an accurate enough solution, the sampling of the semi-dis\-cret\-i\-za\-tion method should be relatively fine $(30 < r <100)$ depending on the actual problem. However, storing all these data points could require large memory. If the $\tau$ deep history is stored into, say, $N_\tau \approx10$ time-steps only then the memory requirements can be reduced. 
	
	\subsection{Initial condition for next simulation}
	As it was defined in Sec.~\ref{sec:iterativeproc}, the selection of proper ICHs has an essential effect on the efficiency of the algorithm.
	
	\begin{figure}
		\centering
		\includegraphics[width=\linewidth]{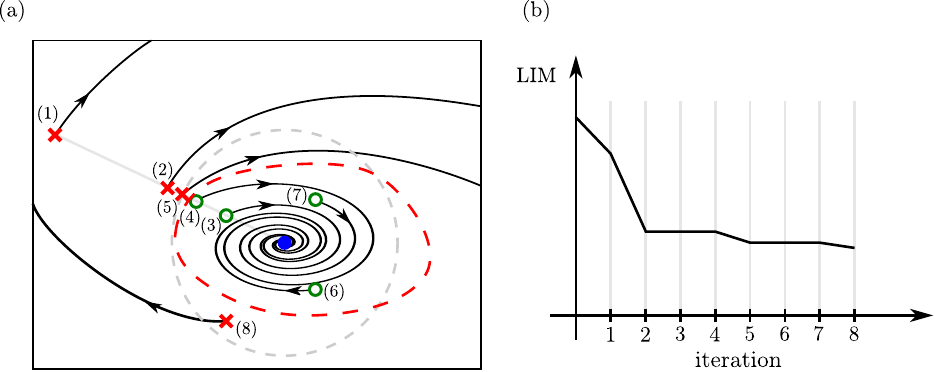}
		\caption{Illustrative example of the algorithm operation. In panel (a), the blue dot represents the desired fixed point, which is encircled by an unstable periodic orbit represented by the red dashed curve. The red crosses and green circles correspond to the ICH of diverging and converging trajectories respectively, while gray dashed circle represents the estimated hypersphere of convergence after 5 steps. Panel (b) presents the evolution of the estimated value of the LIM during the iteration.}
		\label{fig:initcondsfornext}
	\end{figure}
	
	Three methods are suggested for the selection of the ICH of the subsequent simulations: 
	\begin{itemize}
		\item Select a random point within the previously estimated hypersphere of convergence. This is useful to detect new divergent areas and to reduce the hypersphere accordingly.
		\item Apply the bisection method between a divergent ICH and the desired fixed point. This leads to a trajectory starting near the boundary of the BoA, which is very informative.
		\item Choose the next ICH where the last diverging trajectory is closest to the desired fixed point.
	\end{itemize}
	
	These three methods can be combined; in particular, the first few ICHs are randomly chosen until a diverging trajectory is found. Then, it is worth to start a bisection method, which can effectively reduce the estimated hypersphere of convergence. However, the bisection method tends to focus on a specific region, so to study more uniformly the phase space, some randomness must be introduced later. The solutions corresponding to ICHs within the BoA are  converging and do not reduce the estimated LIM; therefore, the probability function for random choice of ICH should have larger probability in the vicinity of the boundary of the estimated hypersphere of convergence, where trajectories are more likely to diverge.
	
	Finally, selecting the ICH at a point of a divergent trajectory which is closest to the desired fixed point, may also be useful. Since, the new time integration starts from a constrained type of initial function, while the corresponding part of the previous trajectory was calculated for the whole nonlinear time-delayed system, the new trajectory will differ from the previous one. Since it is diverging with a large probability, it is a good candidate to reduce the estimated LIM.
	
	Figure~\ref{fig:initcondsfornext}(a) shows an example  of the above strategy, while the iteration towards the corresponding estimated LIM is visualized in panel (b). The blue dot refers to the desired fixed point, while the red dashed curve is an unstable limit cycle, which characterizes the BoA. The first randomly selected headpoint yields a diverging trajectory, after which a bisection method of length four starts (headpoints (2)-(5)). The ICH of the last diverging trajectory determines the 5th estimation of the LIM; the corresponding hypershpere of convergence is visualised as gray dashed circle. After that, the randomly selected headpoints (6) and (7) correspond to converging trajectories, thus they do not reduce the estimated LIM. In contrary, (8) does; thus, another bisection method should be started between that point and the desired fixed point.
	
	
	%
	
	\section{Case studies}\label{sec:casestudies}
	
	In this section, four case studies are presented. First, a Duffing oscillator is considered with delayed linear term and cubic nonlinearity, the basin of attraction of which was presented in Fig.~\ref{fig:DuffingBoA}. Then, the dynamics of a 1 DoF turning operation is investigated, which is later expanded with a nonlinear tuned vibration absorber yielding a 2 DoF system. Finally, the algorithm is applied to the delayed proportional-derivative control of an inverted pendulum.
	
	
	\subsection{Delayed Duffing oscillator}\label{sec:Duffing}
	
	The dimensionless governing equation of the chosen unforced delayed Duffing oscillator is
	\begin{equation}
	\ddot{x}(t)+2\zeta\dot{x}(t)-x(t-\tau)+ax^3(t) = 0\,,
	\end{equation}
	where $\zeta$ is the damping ratio, $a$ is the nonlinear term coefficient and $\tau$ denotes the time delay appearing in the linear term; the LIM of the undelayed system was analysed in \cite{habib2021dynamical}.
	
	The three equilibria of the system are located at $\mathbf{x}_{1}=(x_1,\dot{x}_1)=(-1/\sqrt{a},0)$, $\mathbf{x}_{2}=(x_2,\dot{x}_2)=(0,0)$ and at $\mathbf{x}_{3}=(x_3,\dot{x}_3)=(1/\sqrt{a},0)$, from which $\mathbf{x}_{2}$ is unstable while the other two are stable. We consider that $\mathbf{x}_{1}$ is the desired fixed point, while the other stable fixed point $\mathbf{x}_{3}$ is an undesired one. 
	
	In the preprocessing stage, the location of the desired equilibrium  is provided to the algorithm, while the other stable equilibrium is purposely not provided, to test the algorithm's ability to identify it. The space boundary is set to $-5$ and $5$ for both the dimensionless position and velocity coordinates. The number of iterations  (corresponding to the number of time series simulated) for given parameters is set to $50$, while the maximal length of each simulation is $1000$ time units.
	Next, the spatial and temporal discretization parameters are specified. The number of equally distributed cells in each direction is $n_\mathrm{disc}=501$, while the sampling delay number within the semi-discretization method is $r=30$.
	
	Since the linear undamped, undelayed system is a simple 1 DoF oscillator, the original coordinates are already modal coordinates. The mapping matrix of the semi-discretization is calculated according to Eqs.~\eqref{eq:semidiscmap} and \eqref{eq:semidiscmatr} using the system parameters $a=1$, $\zeta=0.1$ and $\tau=0.1$ . The initial condition is chosen  to be a free undamped vibration, after which the iteration starts.
	
	\begin{figure}
		\centering
		\includegraphics[width=\linewidth]{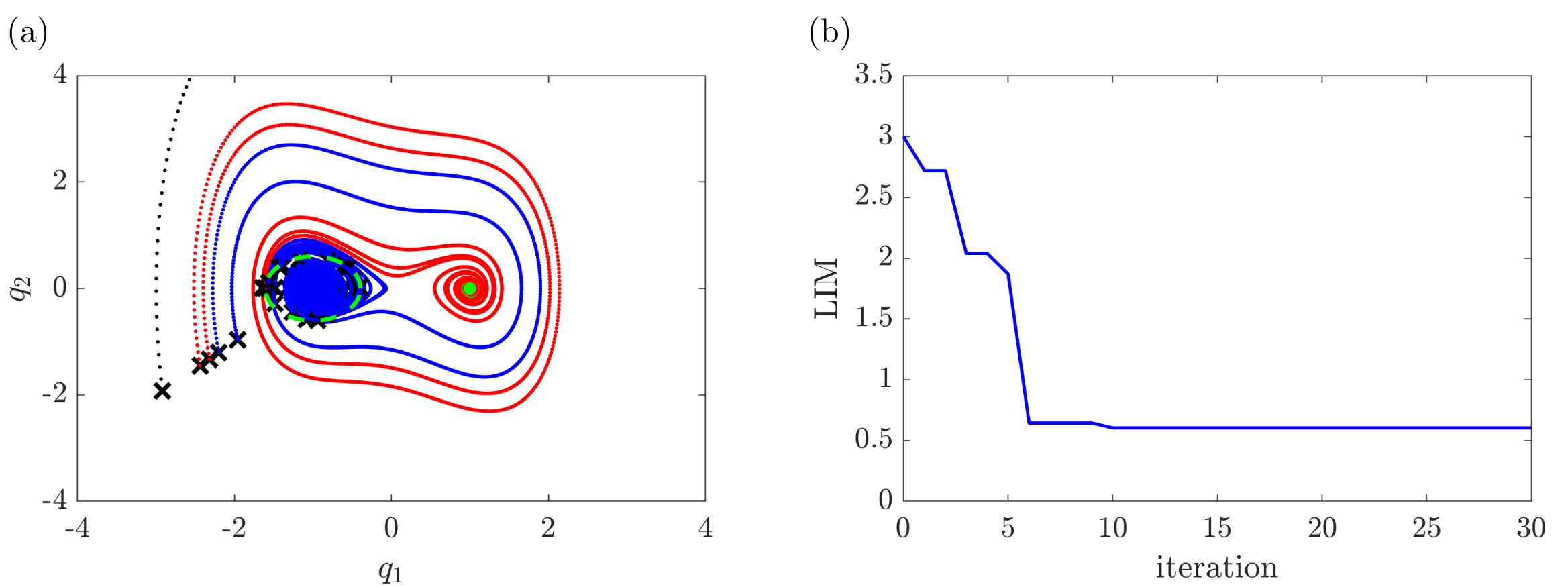}
		\caption{Simulation result in case of the delayed Duffing oscillator $(a=1,~\zeta=0.1,~\tau=0.1)$. In panel (a) blue dots, red dots and black dots refer to trajectories converging to the desired fixed point, converging to the undesired fixed point and diverging out of the space boundary, respectively.}
		\label{fig:Duffing}
	\end{figure}
	
	Figure~\ref{fig:Duffing} shows the result of the iteration. In the left panel, the trajectories are presented in the phase space of the coordinates:  black dots belong to a trajectory that diverges out of the prescribed space boundary; the series of blue dots are convergent trajectories; while red color belongs to trajectories which converge to the undesired fixed point $\mathbf{x}_3$;  black crosses indicate the initial functions' headpoints. Note that close to the equilibria, the series of points that represent the discretized trajectory are so dense that they appear like solid lines. The algorithm also identifies the undesired fixed point, which is indicated by a green dot, while the green dashed ellipse corresponds to the safe region  (the hypersphere of convergence) according to the LIM. 
	
	The first initial condition belongs to the diverging black trajectory; the corresponding ICH and the desired equilibrium $(-1,0)$ serve as a basis for the upcoming bisection method that is applied in 5 steps to obtain a reasonable resolution.
	Note that in case of the Duffing oscillator, the BoAs of the two stable equilibria are layered alternately as one moves away from the origin (see Fig.~\ref{fig:DuffingBoA}). As a consequence, the bisection method does not find the inner boundary of the BoA. Still, a diverging trajectory initiated in the vicinity of the outer BoA boundary remains close to the boundary for a long time. Starting the subsequent simulation from the point of this diverging trajectory, which is located closest to the desired equilibrium, further improves the estimation of the LIM. The bisection method can be applied again based on that point and the desired equilibrium, and the procedure can be repeated till we cannot find a point on a divergent trajectory, which is closer than the previous ones. This way, one point of the hypersphere of convergence can satisfactorily be approximated. Further randomly selected ICHs in the close neighbourhood of the hypersphere help to check the validity of the above estimation.
	
	Panel (b) of Fig.~\ref{fig:Duffing} presents the corresponding estimated values of the LIM during the iteration. It can be observed that the first five steps also reduces the estimated LIM, but the 6th one, which corresponds to the closest point initial condition, reduces the estimated LIM significantly. Subsequent iterations provide only slight reductions of the estimated LIM.
	
	\begin{figure}
		\centering
		\includegraphics[width=\linewidth]{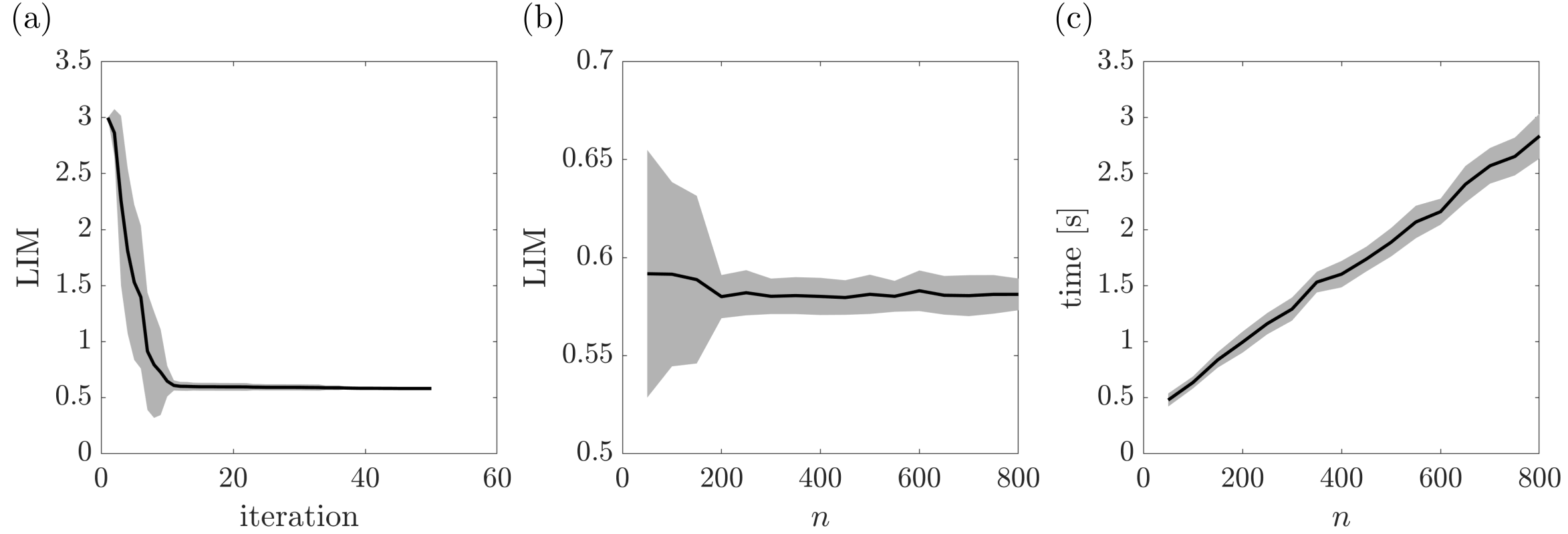}
		\caption{(a) the estimated value of the LIM during an iteration, (b) the estimated LIM varying the spatial discretization number $n$, (c) the required time for a computation for various values of $n$. Each plot is based on 50 independent simulations, the black lines refer to the average, while the gray areas indicate the standard deviation of the results.}
		\label{fig:Duffing_std}
	\end{figure}
	
	Panel (a) of Fig.~\ref{fig:Duffing_std} presents the mean and the standard deviation of the estimated LIM during the iteration calculated from 50 uniformly randomly selected ICHs. One can observe that the iteration converges rapidly, within 10 steps. Panels (b) and (c) show the estimated LIM and the computational time as the spatial discretization number $n$ is increased on the grid. Small values of $n$ yield large standard deviations, while the simulation time increases approximately linearly with $n$. The optimal spatial discretization corresponds to $n\approx200$, where the algorithm provides the estimated LIM in approximately 1 second (processor of the used computer: i7-10510U CPU, 2.30 GHz), which is even faster then the analogous algorithm proposed in \cite{habib2021dynamical} for non-delayed systems.
	
	The change of the sampling delay number $r$ directly modifies the accuracy of the simulations; accordingly, it has an effect on the LIM. However, as the convergence of the approximated dynamics to the real one for increasing $r$ was thoroughly investigated in \cite{insperger2002semi,insperger2011semi}, this aspect is not studied here further.
	

	\subsection{Turning}
	
	\begin{figure}
		\centering
		\includegraphics[width=0.6\linewidth]{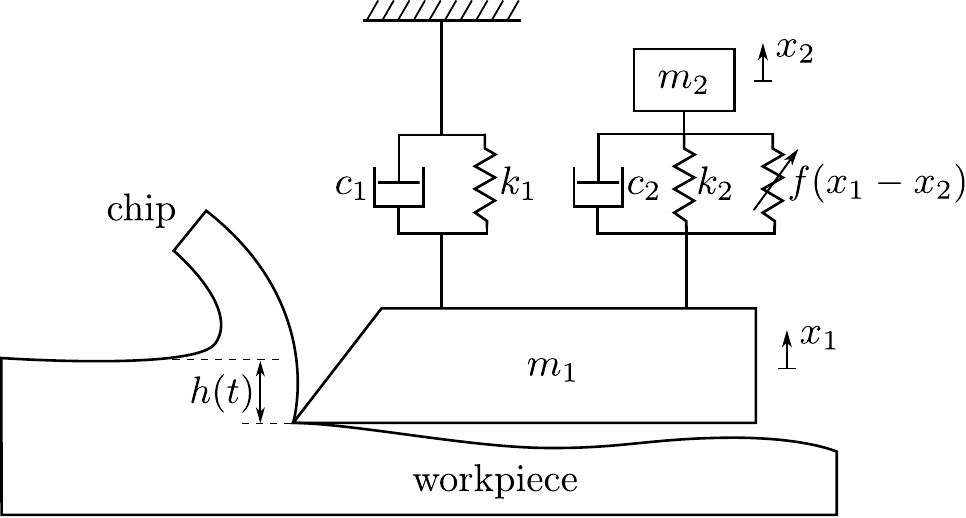}
		\caption{Simplified mechanical model of turning.}
		\label{fig:turningmodel}
	\end{figure}
	
	In the second and third case studies, a turning operation is considered without and with an attached nonlinear tuned vibration absorber (NLTVA), the dynamics of which was investigated in detail in \cite{habib2017chatter,sims2007vibration}. The mechanical model is presented in Fig.~\ref{fig:turningmodel}, the cutting tool of mass $m_1$ is modelled with a single degree of freedom spring-mass-damper system, while the attached NLTVA consists of the mass $m_2$, attached to the cutting tool through a nonlinear spring and a linear damper.  Subscripts $1$ and $2$ refer to the parameters of the cutting tool and to those of the NLTVA, respectively.  The time-dependent depth of cut is denoted by $h(t)$, which takes the value $h_0$ at the equilibrium (steady state cutting). The displacements of the lumped masses $m_1$ and $m_2$ are $x_1$ and $x_2$, respectively, which are measured from the equilibrium position of the system without loss of generality.
	
	First, we consider, the single DoF system without the NLTVA. We introduce the dimensionless time $\tilde{t} = \sqrt{k_1/m_1}\,t$ , and the dimensionless state variables $\tilde{x}_1=x_1/h_0$,  and $\tilde{x}_2=x_2/h_0$. By dropping the tildes, the corresponding dimensionless equation of motion takes the form
	\begin{equation}\label{eq:turning1DoF}
	\begin{split}
	&\ddot{x}_1(t) + 2\zeta_1\dot{x}_1(t) + x_1(t)\\
	&= p((x_1(t-\tau)-x_1(t)) + \eta_2(x_1(t-\tau)-x_1(t))^2+\eta_3(x_1(t-\tau)-x_1(t))^3)\,.
	\end{split}
	\end{equation}
	Here, $\zeta_1=c_1/(2\sqrt{k_1m_1})$ is the damping ratio and the right-hand side of the equation indicates the dimensionless cutting force, which depends on the dimensionless chip width $p$, and on the experimentally determined coefficients $\eta_2$ and $\eta_3$ (see \cite{habib2017chatter} for more detail). Assuming that the tool does not leave the workpiece, that is, no fly-over effect occurs, the actual depth of cut is given as the difference of the position of the cutting tool at the present time and one revolution earlier, which causes the appearance of the delay $\tau$  in the equation. The corresponding dimensionless spindle speed is $\Omega_\mathrm{d}=2\pi/\tau$.  
	
	\begin{figure}
		\centering
		\includegraphics[width=\linewidth]{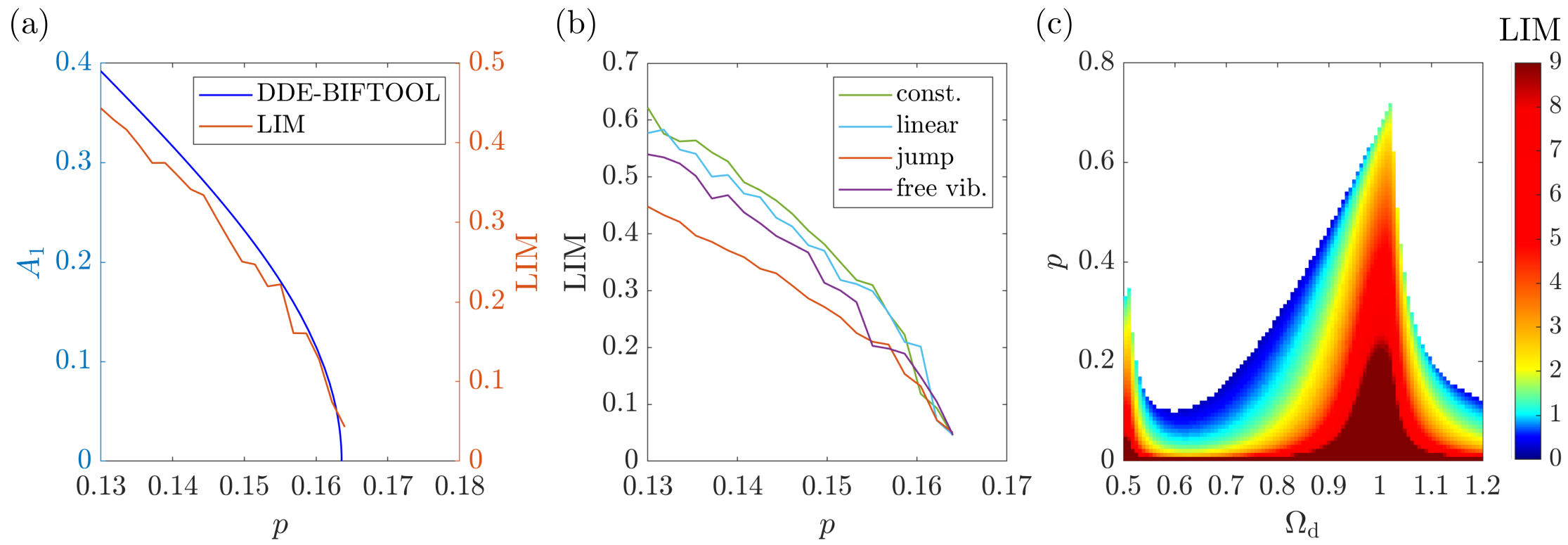}
		\caption{Considering the single DoF model of turning ($\zeta_1=0.05,~\eta_2=-0.5209,~\eta_3=0.6547$), (a) comparison of the bifurcation diagram created with DDE-BIFTOOL and the estimated LIM with free vibration initial condition ($\tau=9$, and so $\Omega_\mathrm{d}=0.6987$), (b) comparison of the four types of initial conditions ($\Omega_\mathrm{d}=0.6987$), (c) stability chart colored according to the estimated value of the LIM. }
		\label{fig:turning1_plim}
	\end{figure}

	The proposed algorithm was applied to estimate the LIM of the desired steady state cutting. First, the damping ratio and the time delay are fixed to $\zeta_1 = 0.05$ and $\tau=9$ and the LIM is calculated for varying dimensionless chip width $p$. For these parameter values, the system loses stability through a subcritical Andronov-Hopf bifurcation. Therefore, an unstable limit cycle exists around the trivial stable solution \cite{dombovari2008estimates}.
	
	In panel (a) of Fig.~\ref{fig:turning1_plim}, the tendency of the estimated LIM is compared to the bifurcation diagram built with DDE-BIFTOOL \cite{engelborghs2002numerical}. Of course, the two curves correspond to different measures, thus, they are quantitatively different, but the algorithm accurately reproduces the trend of the bifurcation diagram. 
	
	Panel (b) of Fig.~\ref{fig:turning1_plim} presents the estimated LIM values in case of different types of initial conditions. From an engineering point of view, they are practically equivalent as they do not affect the general trend of the estimated LIM.
	
	In panel (c) of Fig.~\ref{fig:turning1_plim}, the LIM was calculated sweeping both the dimensionless spindle speed $\Omega_\mathrm{d}$ and the dimensionless chip width $p$ for parameter values where the trivial solution is stable. The figure illustrates that the LIM takes smaller values on the left-hand side of the lobe than on the right-hand side, which means that stationary cutting with parameters at the left side of the lobe is less robust against external perturbations. This is a relevant information from engineering viewpoint, which would be computationally demanding to obtain with alternative methods. 
	
	
	\subsection{Turning with NLTVA}
	In this section, the same turning operation is considered as before, but now the effect of the NLTVA is also taken into account (see Fig.~\ref{fig:turning1_plim}). Note that the enlargement of the system dimension from 2 to 4 can be a significant challenge  \cite{belardinelli2016first}.
	The corresponding dimensionless governing equations take the form:
	\begin{align}
	\begin{split} \label{eq:turning2DoF1}
	&\ddot{x}_1(t) + 2\zeta_1\dot{x}_1(t) +  x_1(t) +2\zeta_2\gamma\mu(\dot{x}_1(t)-\dot{x}_2(t))+\gamma^2\mu(x_1(t)-x_2(t))\\ &+\alpha_3(x_1(t)-x_2(t))^3
	= p((x_1(t-\tau)-x_1(t)) + \eta_2(x_1(t-\tau)-x_1(t))^2 \\ &+\eta_3(x_1(t-\tau)-x_1(t))^3)\,,
	\end{split}\\ \label{eq:turning2DoF2}
	\mu& \ddot{x}_2(t)+ 2\zeta_2\gamma\mu(\dot{x}_2(t)-\dot{x}_1(t)) + \gamma^2\mu(x_2(t)-x_1(t))+\alpha_3(x_2(t)-x_1(t))^3 = 0 ,
	\end{align}
	where $\alpha_3$ is the absorber dimensionless cubic stiffness parameter, while the further dimensionless parameters
	\begin{equation}
	\zeta_2 = \frac{c_2}{2\sqrt{k_2m_2}}, \quad \mu=\frac{m_2}{m_1}\,,\quad \mathrm{and} \quad \gamma = \frac{\omega_{\mathrm{n},2}}{\omega_{\mathrm{n},1}} = \sqrt{\frac{k_1m_2}{m_1k_2}}\,, 
	\end{equation}
	are the damping ratio of the NLTVA, the mass ratio, and the frequency ratio, respectively. As before, the right-hand side of Eq.~\eqref{eq:turning2DoF1} presents the dimensionless cutting force that contains also delayed terms.
	
	\begin{figure}
		\centering
		\makebox[\textwidth][c]{\includegraphics[width=1\textwidth]{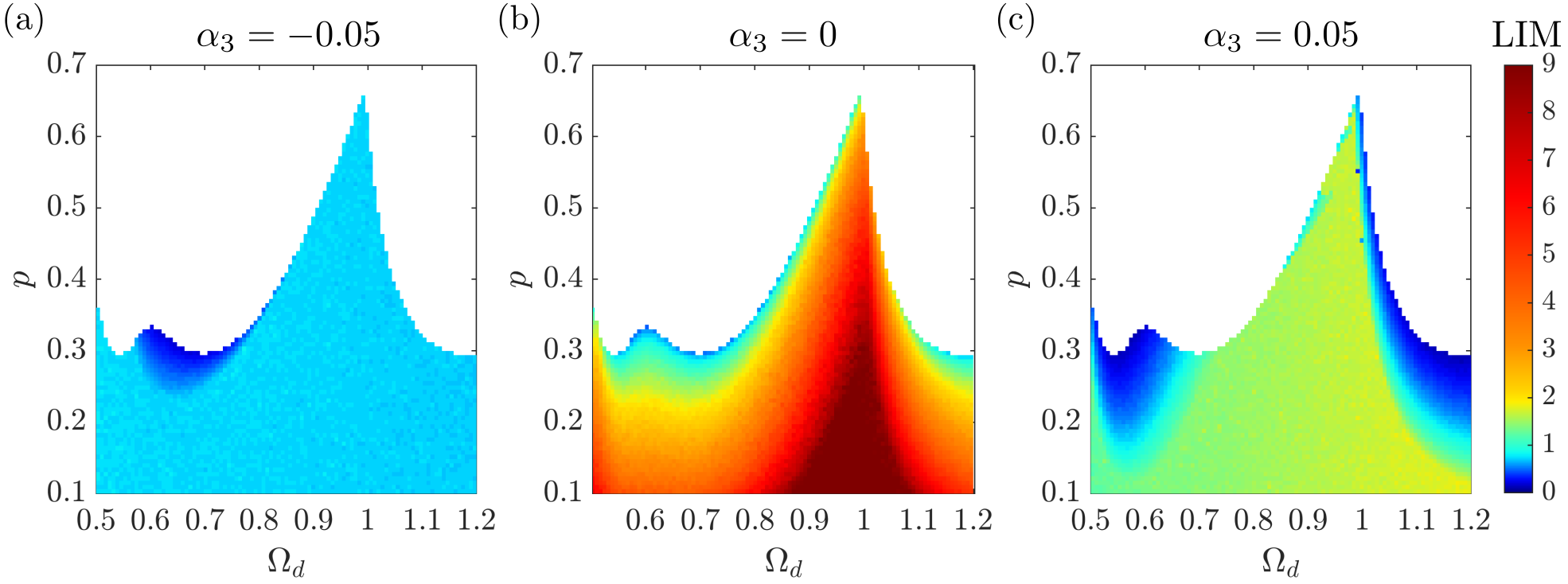}}%
		\caption{Stability charts of the 2 DoF turning operation for various values of the absorber cubic stiffness parameter $\alpha_3$ with the estimated LIM indicated by the color scheme for parameters $\zeta_1=0.05,~\mu=0.05,~\gamma=1.069,~\zeta_2=0.1437,~\eta_2=-0.5209,~\eta_3=0.6547$. }
		\label{fig:turning2omdlim}
	\end{figure}
	
	The nonlinear spring of the NLTVA is described with a softening (hardening) characteristic as the value of the cubic stiffness parameter $\alpha_3$ is negative (positive). After fixing the parameter values as indicated in the caption of Fig.~\ref{fig:turning2omdlim}, the corresponding dynamics were analysed with the proposed algorithm, and the LIM was calculated in the plane of the dimensionless spindle speed $\Omega_\mathrm{d}$ and chip width $p$ for various values of the absorber cubic stiffness parameter $\alpha_3$ (see Fig.~\ref{fig:turning2omdlim}). Panel (b) presents the case of a linear absorber, which provides a robust system for most of the parameter combinations; however, in the vicinity of the stability boundary, the dynamical integrity is still limited due to the subcritical nature of the Andronov-Hopf bifurcation occurring there.
	In \cite{habib2017chatter}, it was proven that including a nonlinear term in the stiffness function, either softening or hardening, can change the character of the bifurcation to a supercritical one in some parts of the stability boundary.
	Panels (a) and (b) of Fig.~\ref{fig:turning2omdlim} illustrate the dynamical integrity for a softening and hardening spring, respectively.
	We note that the overall dynamical integrity is reduced in both cases, especially in the softening case, regardless of the bifurcations' criticality \cite{habib2017chatter}. 
	This observation suggests that the vibration absorber should be kept as linear as possible and that a local bifurcation analysis can be misleading concerning global dynamics. 
	
	In this four dimensional time delayed case, the algorithm can calculate the dynamical integrity  corresponding to 100 points of the parameter space within 4 minutes using a standard laptop without parallel computation (i7-10510U CPU, 2.30 GHz). This allows the algorithm to be used for parametric analysis, while these results are practically infeasible to be produced with other methods. 
	
	In Fig.~\ref{fig:turning2omdlim}(c), one can observe a light blue area, that is, an area with smaller LIM, along the stability boundary at around $\Omega_\mathrm{d}=0.9$, $p=0.5$. In these points, the system has a very small linear damping, and if the ICH was close to the boundary of the hypersphere of convergence, then the corresponding trajectories converged so slowly to the equilibrium, that they did not reach it within the predefined time interval. Thus, the algorithm categorised them as diverging trajectories, yielding erroneously smaller values of the LIM than the correct ones obtained by extended simulations for longer periods. 
	
	\subsection{Inverted pendulum with NLTVA}\label{sec:invpend}
	\begin{figure}
		\centering
		\includegraphics[width=\linewidth]{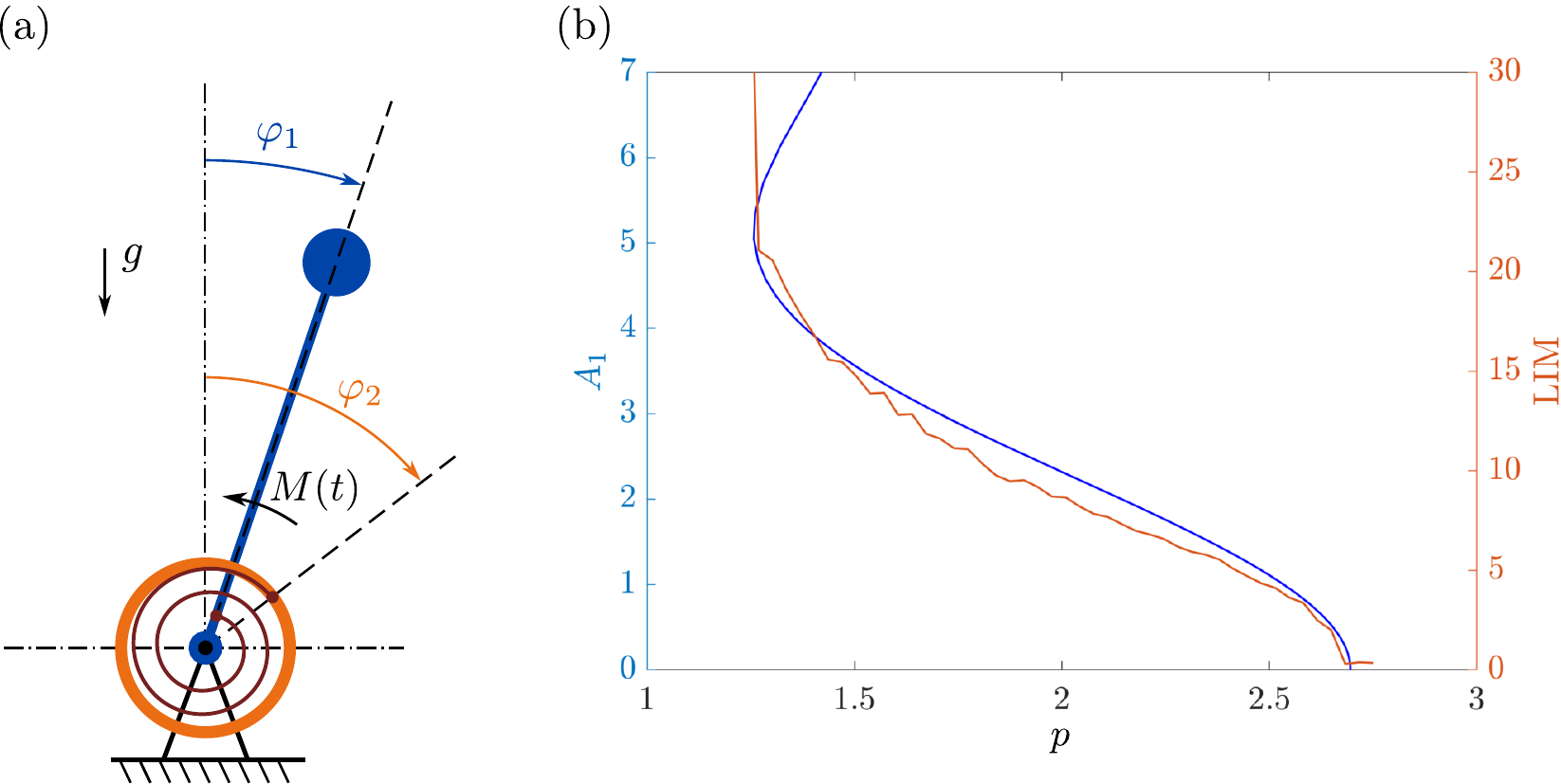}
		\caption{(a) Model of the inverted pendulum with NLTVA. (b) Comparison of the bifurcation diagram created with DDE-BIFTOOL and the estimated LIM ($\mu=0.1$, $\gamma=2.3$,  $\zeta_2=0.174$, $d=2.8$, $\tau=0.5$).}
		\label{fig:pendplim}
	\end{figure}
	
	The fourth case study is the analysis of the dynamics of a controlled inverted pendulum with an attached NLTVA in the gravitational field. The mechanical model is presented in the left panel of Fig.~\ref{fig:pendplim}. The positions of the pendulum and the NLTVA are described by the angles $\varphi_1$ and $\varphi_2$, respectively. A delayed proportional-derivative control torque $M(t)=P\varphi_1(t-\tau)+D\dot{\varphi}_1(t-\tau)$ is applied to the pendulum to keep it at the upright position.
	The corresponding dimensionless governing equations assume the form
	\begin{align}    
	\ddot{\varphi}_1-\sin(\varphi_1)+2\zeta_2\mu\gamma(\dot{\varphi}_1-\dot{\varphi}_2) + \mu\gamma^2(\varphi_1-\varphi_2) 
	&= - p\varphi_1(t-\tau)-d\dot{\varphi}_1(t-\tau) ,
	\\
	\mu(\ddot{\varphi}_2 + 2\zeta_2\gamma(\dot{\varphi}_2-\dot{\varphi}_1) + \gamma^2(\varphi_2-\varphi_1)) &= 0,
	\end{align}
	where $\zeta_2$ is the relative damping of the NLTVA, $\mu$ is the ratio between the mass moment of inertia of the NLTVA and that of the inverted pendulum, $\gamma$ denotes the natural frequency ratio, while $p$ and $d$ are the dimensionless proportional and differential gains.
	
	The corresponding nonlinear dynamics was analysed in \cite{habib2019suppression} in detail proving that the Andronov-Hopf bifurcation at the stability boundary is always subcritical. 
	In Fig.~\ref{fig:pendplim}, the result of the proposed algorithm (red curve) is compared to the bifurcation diagram created with DDE-BIFTOOL (blue curve) \cite{engelborghs2002numerical}. Choosing the dimensionless proportional gain as the bifurcation parameter and fixing other parameters at $\mu=0.1$, $\gamma=2.3$, $\zeta_2=0.174$, $d=2.8$, and $\tau=0.5$, the linear system is stable for $p<2.695$. The branch of periodic solutions bends above the linearly stable region, in accordance with the subcriticality of the bifurcation, yielding an unstable limit cycle around the stable trivial fixed point. Later, the branch of the limit cycle has a fold point at $p=1.257$, where the curve bends back to the right and the related periodic solutions become stable. This implies the existence of a stable limit cycle coexisting with an unstable one and the stable equilibrium. 
	
	Again, the LIM was determined for free vibration initial conditions. However, the weights of the distance definition in Eq.~\eqref{eq:distance} could not be calculated with the energy based method proposed in \eqref{eq:distancedef} since the fixed point is unstable when the delayed terms are neglected. One way to solve this problem is to keep the proportional control term without time delay, and calculate the distance based on that. However, it yields that the distance measure is a function of $p$, while $p$ is the varying bifurcation parameter. Accordingly, for this case, the vector of weights (see Eq.~\ref{eq:distanceweights}) was simply chosen to be $\alpha=[1,1,1,1]$.
	
	As it can be observed in Fig.~\ref{fig:pendplim}, the shape of the LIM curve matches well with the bifurcation diagram. Moreover, it accurately predicts the fold point, below which the trivial fixed point is globally stable, that is, the LIM is theoretically infinite, while in practice, it is limited by the predefined spatial boundary of the algorithm. 
	
\begin{figure}[t!]
		\centering    \includegraphics[width=\linewidth]{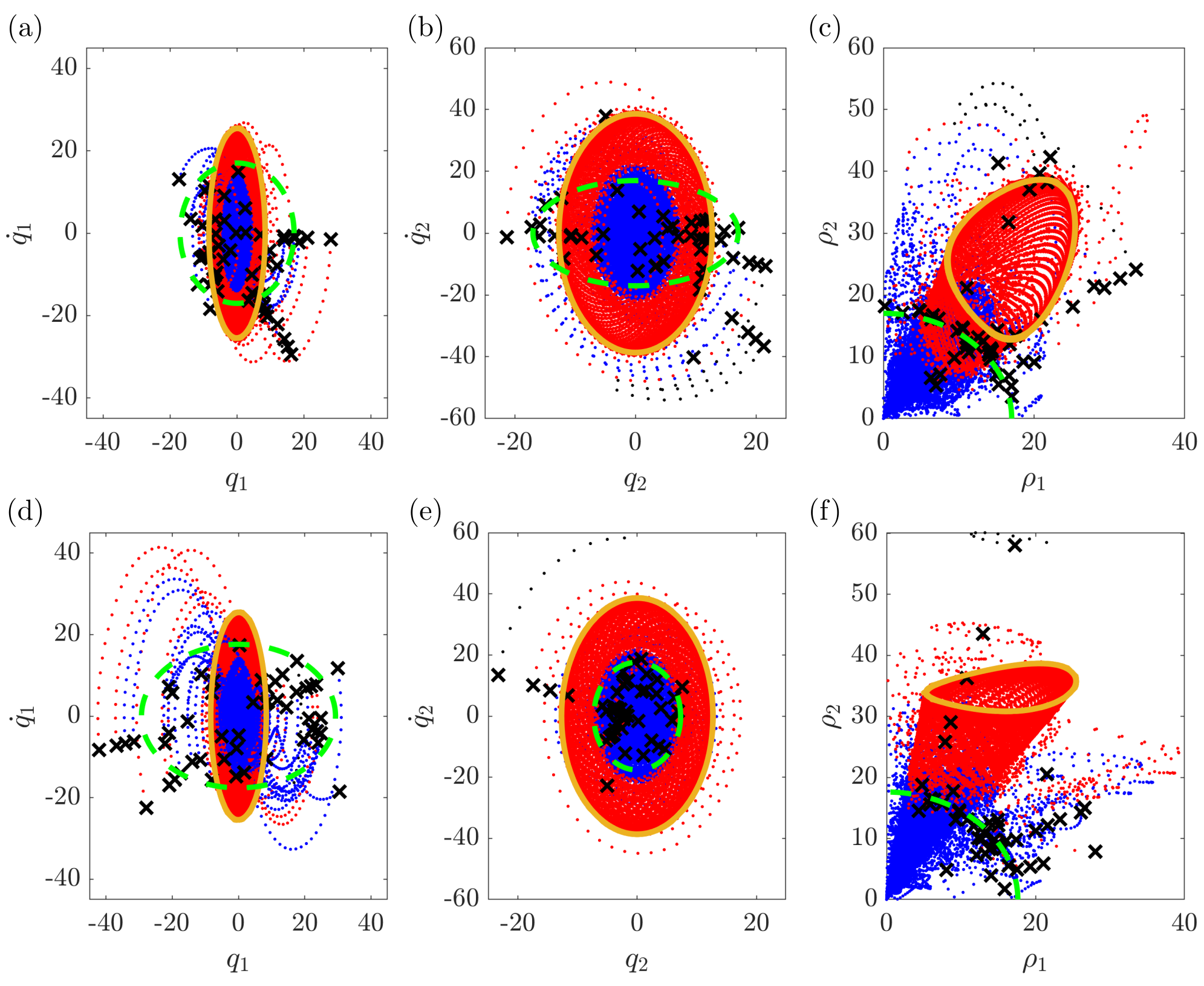}
		\caption{Panels (a,b,d,e) present the projection of the trajectories onto the phase space of the inverted pendulum, while panels (c,f) show the trajectories in the plane of the square root of modal energies ($\mu=0.1$, $\gamma=2.3$,  $\zeta_2=0.174$, $d=2.8$, $p=1.4$, $\tau=0.5$). Panels (a,b,c) present simulations with  user defined distance weights $\alpha=[1,1,1,1]$; (d,e,f) present simulations when the weights were calculated taking into account the delay-free proportional control term. Blue dots: converging trajectories, red dots: non-converging trajectories, orange curve: projection of the stable limit cycle, green dashed curve: projection of the hypersphere of convergence.}
		\label{fig:pendq1q2}
	\end{figure}

	Figure~\ref{fig:pendq1q2} presents the simulated trajectories of the algorithm projected onto the two dimensional subspaces of the system dynamics. Here, $q_i$ and $\dot{q}_i$ ($i=1,2$) are the modal coordinates and their derivatives, while the variables  $\rho_i = \sqrt{\omega_i^2q_i^2+\dot{q}_i^2}$ ($i=1,2$) are related to the energy of the corresponding mode. Note that the distance definition in Eq.~\eqref{eq:distancedef} can be reformulated as $d=\sqrt{\rho_1^2+\rho_2^2}$.
	Panels (a,b,c) correspond to the user defined weight vector $\alpha=[1,1,1,1]$, while panels (d,e,f) refer to the weights calculated with the inclusion of the proportional control term. Blue dots denote the convergent while red dots the divergent trajectories. The projection of the hypersphere of convergence is the green dashed curve, while the projection of the stable limit cycle is presented by the orange curve.
	
	The good qualitative matching between the bifurcation diagram and the LIM trend in Fig.~\ref{fig:pendplim} suggests that the unit weights do not effect the trend of the LIM in the parameter space in general. However, Fig.~\ref{fig:pendq1q2} illustrates that the weights do affect the absolute value of the LIM. In particular, the velocity is more important than the displacement in the given example. 
	
	Selecting the weights based on the natural angular frequencies, the hypersphere of convergence typically fits the trajectories much better and gives equal weights to the positions and velocities. However, in this case, this is valid only to the second mode, because the stable limit cycle has a vibration frequency close to the second natural frequency. As panel (f) shows, the second mode carries most of the energy of the system. Nonetheless, it should be remembered that the weights are based on linear properties including the non-delayed proporional gain, while the trajectories are significantly affected by the non-linearities. 
	
	We note that the projections in Fig.~\ref{fig:pendq1q2} might give the impression that many divergent trajectories cross the hypersphere of convergence. However, this is simply due to the projection of a 4-dimensional space onto a plane. In reality, diverging trajectories enter the hypersphere of convergence only rarely and slightly because of the definition of BoA based on the ICHs discussed above.
	
	\section{Discussion and conclusions}\label{sec:discussion}
	The case studies show the capabilities of the proposed algorithm, which gives a rapid and accurate estimation of the dynamical integrity of a stable fixed point. Moreover, the algorithm enables to compute the LIM so fast that it can be used for analysing the variation of the dynamical integrity in the parameter space. Other methods are far more expensive and require the use of supercomputers.
	
	The developed algorithm is based on the following concepts:  
	\begin{itemize}
		\item The LIM is a quantitative measure that characterizes the robustness of the examined fixed point. By definition, the corresponding hypersphere of convergence does not include fractal-like domains and it gives a conservative  measure of the safe region.
		\item The semi-discretization provides a fast and accurate scheme for determining the trajectories of nonlinear delayed dynamical systems.
	\end{itemize}
	
	The main advantages of the algorithm are as follows:
	\begin{itemize}
		\item The algorithm is generalized for applications in time delayed systems. 
		\item The estimated value of the LIM significantly reduces  in the first few steps, enabling the algorithm to converge rapidly to a reasonable estimate of the LIM.
		\item The algorithm detects whether a trajectory converges to an already categorized one, reducing the length of new simulations with the number of iterations.
		\item There are no memory issues, as it is not required to investigate all the cells of the phase space, not even all the cells within the hypersphere of convergence. Conversely, the cell mapping method has significant memory issues.
	\end{itemize}
	Despite these advantages, the proposed algorithm has some limitations:
	\begin{itemize}
		\item Although, the algorithm can handle any kind of discretized initial functions, the estimated LIM is restricted to a fixed type of initial function.  The case studies show that the type of initial function changes the absolute value of LIM, but its trend in the parameter space is qualitatively unaffected.
		\item The estimated LIM is reduced if a diverging trajectory is found within the estimated hypersphere of convergence, meaning that the algorithm gives an upper estimate of the LIM. Thus, it gives a non-conservative estimate of the conservative measure of the safe robust region.
		\item Large dimensional phase spaces are filled relatively slowly, as trajectories are single dimensional objects, regardless of the phase space dimension. Therefore, undesired stable solutions within the hypersphere of convergence may become undetected, unless the number of iterations is significantly increased for large dimensional systems. This limitation can be mitigated resorting to model reduction techniques \cite{scarciotti2017nonlinear,haller2016nonlinear,chaturantabut2010nonlinear}.
	\end{itemize}
	
	Four case studies were presented in which the advantages and limitations of the proposed algorithm are discussed. For each mechanical system, the algorithm provided a fast and accurate estimation of LIM. In case of parametric analysis, the results showed good agreement in the qualitative trends with the bifurcation diagrams. The relative rapidity of the algorithm enables researchers and engineers to use it for parametric analysis on a sufficiently large scale, making it applicable also for design.
	
	In the future, the algorithm is planned to go through further developments. Currently, if the algorithm is used for parametric analysis, then the LIM is determined separately for each value of the varying system parameter. However, in most cases, the variation of the LIM is smooth, or at least piece-wise continuous. Thus, the previous estimation of the LIM could be used as an initial guess for the next iteration, which  speeds up the algorithm.
	
	The proposed algorithm currently investigates the dynamical integrity of stable fixed points. However, we plan to extend it to other types of steady-state solutions, such as periodic ones, which does not require significant changes to the main structure of the algorithm.

	\section*{Code availability}
	The code supporting the results presented in the paper is publicly available at \url{https://gdrg.mm.bme.hu}
	
	\section*{Acknowledgements}
	The research reported in this paper has been supported by the National Research, Development and Innovation Office (Grant No. NKFI-134496, NKFI-KKP-133846 and NKFI-K-132477).  
	
	
	
	\bibliographystyle{elsarticle-num} 

\begin{thebibliography}{10}
\expandafter\ifx\csname url\endcsname\relax
  \def\url#1{\texttt{#1}}\fi
\expandafter\ifx\csname urlprefix\endcsname\relax\def\urlprefix{URL }\fi
\expandafter\ifx\csname href\endcsname\relax
  \def\href#1#2{#2} \def\path#1{#1}\fi

\bibitem{pacejka1965analysis}
H.~B. Pacejka, Analysis of the shimmy phenomenon, Proceedings of the
  Institution of Mechanical Engineers: Automobile Division 180~(1) (1965)
  251--268.

\bibitem{beregi2019bifurcation}
S.~Beregi, D.~Takacs, G.~Stepan, Bifurcation analysis of wheel shimmy with
  non-smooth effects and time delay in the tyre--ground contact, Nonlinear
  Dynamics 98~(1) (2019) 841--858.

\bibitem{horvath2022stability}
H.~Z. Horvath, D.~Takacs, Stability and local bifurcation analyses of
  two-wheeled trailers considering the nonlinear coupling between lateral and
  vertical motions, Nonlinear Dynamics (2022) 1--18.

\bibitem{habib2023towed}
G.~Habib, A.~Epasto, Towed wheel shimmy suppression through a nonlinear tuned
  vibration absorber, Nonlinear Dynamics (2023) 1--14.

\bibitem{lind2003flight}
R.~Lind, Flight-test evaluation of flutter prediction methods, Journal of
  Aircraft 40~(5) (2003) 964--970.

\bibitem{verstraelen2017experimental}
E.~Verstraelen, G.~Habib, G.~Kerschen, G.~Dimitriadis, Experimental passive
  flutter suppression using a linear tuned vibration absorber, AIAA Journal
  55~(5) (2017) 1707--1722.

\bibitem{basta2021flutter}
E.~Basta, M.~Ghommem, S.~Emam, Flutter control and mitigation of limit cycle
  oscillations in aircraft wings using distributed vibration absorbers,
  Nonlinear Dynamics 106 (2021) 1975--2003.

\bibitem{takarics2020active}
B.~Takarics, B.~Patartics, T.~Luspay, B.~Vanek, C.~Roessler, J.~Bartasevicius,
  S.~J. Koeberle, M.~Hornung, D.~Teubl, M.~Pusch, et~al., Active flutter
  mitigation testing on the flexop demonstrator aircraft, in: AIAA Scitech 2020
  Forum, 2020, p. 1970.

\bibitem{drachinsky2022flutter}
A.~Drachinsky, O.~Avin, D.~E. Raveh, Y.~Ben-Shmuel, M.~Tur, Flutter tests of
  the pazy wing, AIAA Journal 60~(9) (2022) 5414--5421.

\bibitem{kerswell2018nonlinear}
R.~Kerswell, Nonlinear nonmodal stability theory, Annual Review of Fluid
  Mechanics 50 (2018) 319--345.

\bibitem{cherubini2015nonlinear}
S.~Cherubini, P.~De~Palma, J.-C. Robinet, Nonlinear optimals in the asymptotic
  suction boundary layer: Transition thresholds and symmetry breaking, Physics
  of Fluids 27~(3) (2015) 034108.

\bibitem{nagatani2002physics}
T.~Nagatani, The physics of traffic jams, Reports on progress in physics 65~(9)
  (2002) 1331.

\bibitem{orosz2010traffic}
G.~Orosz, R.~E. Wilson, G.~St{\'e}p{\'a}n, Traffic jams: dynamics and control
  (2010).

\bibitem{molnar2021delayed}
T.~G. Moln{\'a}r, D.~Upadhyay, M.~Hopka, M.~Van~Nieuwstadt, G.~Orosz, Delayed
  lagrangian continuum models for on-board traffic prediction, Transportation
  Research Part C: Emerging Technologies 123 (2021) 102991.

\bibitem{veraszto2017nonlinear}
Z.~Veraszto, G.~Stepan, Nonlinear dynamics of hardware-in-the-loop experiments
  on stick--slip phenomena, International Journal of Non-Linear Mechanics 94
  (2017) 380--391.

\bibitem{szaksz2022transient}
B.~Szaksz, G.~Stepan, Transient chaotic behavior of fuzzy controlled polishing
  processes, Chaos: An Interdisciplinary Journal of Nonlinear Science 32~(9)
  (2022) 093112.

\bibitem{habib2022bistability}
G.~Habib, A.~B{\'a}rtfai, A.~Barrios, Z.~Dombovari, Bistability and delayed
  acceleration feedback control analytical study of collocated and
  non-collocated cases, Nonlinear Dynamics 108~(3) (2022) 2075--2096.

\bibitem{bartfai2022hopf}
A.~Bartfai, Z.~Dombovari, Hopf bifurcation calculation in neutral delay
  differential equations: Nonlinear robotic arms subject to delayed
  acceleration feedback control, International Journal of Non-Linear Mechanics
  147 (2022) 104239.

\bibitem{dombovari2015bistable}
Z.~Dombovari, G.~Stepan, On the bistable zone of milling processes,
  Philosophical Transactions of the Royal Society A: Mathematical, Physical and
  Engineering Sciences 373~(2051) (2015) 20140409.

\bibitem{molnar2019closed}
T.~G. Molnar, T.~Insperger, G.~Stepan, Closed-form estimations of the bistable
  region in metal cutting via the method of averaging, International Journal of
  Non-Linear Mechanics 112 (2019) 49--56.

\bibitem{iklodi2022bi}
Z.~Iklodi, D.~A. Barton, Z.~Dombovari, Bi-stability induced by motion limiting
  constraints on boring bar tuned mass dampers, Journal of Sound and Vibration
  517 (2022) 116538.

\bibitem{pourbeik2006anatomy}
P.~Pourbeik, P.~S. Kundur, C.~W. Taylor, The anatomy of a power grid
  blackout-root causes and dynamics of recent major blackouts, IEEE Power and
  Energy Magazine 4~(5) (2006) 22--29.

\bibitem{gajduk2014stability}
A.~Gajduk, M.~Todorovski, L.~Kocarev, Stability of power grids: An overview,
  The European Physical Journal Special Topics 223~(12) (2014) 2387--2409.

\bibitem{smith2017basins}
V.~A. Smith, T.~E. Lockhart, M.~L. Spano, Basins of attraction in human
  balance, The European Physical Journal Special Topics 226 (2017) 3315--3324.

\bibitem{zakynthinaki2010modeling}
M.~S. Zakynthinaki, J.~R. Stirling, C.~A. Cordente~Mart{\'\i}nez, A.~L.
  D{\'\i}az~de Durana, M.~S. Quintana, G.~R. Romo, J.~S. Molinuevo, Modeling
  the basin of attraction as a two-dimensional manifold from experimental data:
  Applications to balance in humans, Chaos: An Interdisciplinary Journal of
  Nonlinear Science 20~(1) (2010) 013119.

\bibitem{aguirre2014bifurcations}
P.~Aguirre, J.~D. Flores, E.~Gonz{\'a}lez-Olivares, Bifurcations and global
  dynamics in a predator--prey model with a strong allee effect on the prey,
  and a ratio-dependent functional response, Nonlinear Analysis: Real World
  Applications 16 (2014) 235--249.

\bibitem{arancibia2019basins}
C.~Arancibia-Ibarra, The basins of attraction in a modified
  may--holling--tanner predator--prey model with allee affect, Nonlinear
  Analysis 185 (2019) 15--28.

\bibitem{rega2021global}
G.~Rega, V.~Settimi, Global dynamics perspective on macro-to nano-mechanics,
  Nonlinear Dynamics 103~(2) (2021) 1259--1303.

\bibitem{martiniani2016structural}
S.~Martiniani, K.~J. Schrenk, J.~D. Stevenson, D.~J. Wales, D.~Frenkel,
  Structural analysis of high-dimensional basins of attraction, Physical Review
  E 94~(3) (2016) 031301.

\bibitem{sprott2015classifying}
J.~Sprott, A.~Xiong, Classifying and quantifying basins of attraction, Chaos:
  An Interdisciplinary Journal of Nonlinear Science 25~(8) (2015) 083101.

\bibitem{yan2021statistical}
Y.~Yan, J.~Xu, M.~Wiercigroch, Q.~Guo, Statistical basin of attraction in
  time-delayed cutting dynamics: Modelling and computation, Physica D:
  Nonlinear Phenomena 416 (2021) 132779.

\bibitem{lenci2019global}
S.~Lenci, G.~Rega, Global Nonlinear Dynamics for Engineering Design and System
  Safety, Vol. 588, Springer, 2019.

\bibitem{ratschan2010providing}
S.~Ratschan, Z.~She, Providing a basin of attraction to a target region of
  polynomial systems by computation of lyapunov-like functions, SIAM Journal on
  Control and Optimization 48~(7) (2010) 4377--4394.

\bibitem{grinberg2015boundary}
I.~Grinberg, O.~V. Gendelman, Boundary for complete set of attractors for
  forced--damped essentially nonlinear systems, Journal of Applied Mechanics
  82~(5) (2015) 051004.

\bibitem{biemond2014estimation}
J.~B. Biemond, W.~Michiels, Estimation of basins of attraction for controlled
  systems with input saturation and time-delays, IFAC Proceedings Volumes
  47~(3) (2014) 11006--11011.

\bibitem{hsu1980theory}
C.~S. Hsu, A theory of cell-to-cell mapping dynamical systems, Journal of
  Applied Mechanics 47~(4) (1980) 931--939.

\bibitem{hsu1986cell}
C.~S. Hsu, H.~M. Chiu, {A Cell Mapping Method for Nonlinear Deterministic and
  Stochastic Systems—Part I: The Method of Analysis}, Journal of Applied
  Mechanics 53~(3) (1986) 695--701.

\bibitem{hsu2013cell}
C.~S. Hsu, Cell-to-cell mapping: a method of global analysis for nonlinear
  systems, Vol.~64, Springer Science \& Business Media, 2013.

\bibitem{sun2018cell}
J.-Q. Sun, F.-R. Xiong, O.~Sch{\"u}tze, C.~Hern{\'a}ndez, Cell mapping methods,
  Springer, 2018.

\bibitem{liu2016global}
X.~Liu, L.~Hong, J.~Jiang, Global bifurcations in fractional-order chaotic
  systems with an extended generalized cell mapping method, Chaos: An
  Interdisciplinary Journal of Nonlinear Science 26~(8) (2016) 084304.

\bibitem{andonovski2020six}
N.~Andonovski, S.~Lenci, Six-dimensional basins of attraction computation on
  small clusters with semi-parallelized scm method, International Journal of
  Dynamics and Control 8 (2020) 436--447.

\bibitem{habib2021dynamical}
G.~Habib, Dynamical integrity assessment of stable equilibria: a new rapid
  iterative procedure, Nonlinear Dynamics 106~(3) (2021) 2073--2096.

\bibitem{hu2003dynamics}
H.~Hu, Z.~Wang, D.~Schaechter, Dynamics of controlled mechanical systems with
  delayed feedback, Appl. Mech. Rev. 56~(3) (2003) B37--B37.

\bibitem{szaksz2022delay}
B.~Szaksz, G.~Stepan, Delay-induced bifurcations in collocated position control
  of an elastic arm, Nonlinear Dynamics 107~(2) (2022) 1611--1622.

\bibitem{szaksz2022nonlinear}
B.~Szaksz, G.~Stepan, Nonlinear oscillations in delayed collocated control of
  pendulum on trolley, in: International Design Engineering Technical
  Conferences and Computers and Information in Engineering Conference, Vol.
  86304, American Society of Mechanical Engineers, 2022, p. V009T09A032.

\bibitem{sun2017effect}
X.~Sun, J.~Xu, J.~Fu, The effect and design of time delay in feedback control
  for a nonlinear isolation system, Mechanical Systems and Signal Processing 87
  (2017) 206--217.

\bibitem{de2020memory}
R.~De~Luca, F.~Romeo, Memory effects and self-excited oscillations in
  deterministic epidemic models with intrinsic time delays, The European
  Physical Journal Plus 135~(10) (2020) 1--17.

\bibitem{stepan2001modelling}
G.~St{\'e}p{\'a}n, Modelling nonlinear regenerative effects in metal cutting,
  Philosophical Transactions of the Royal Society of London. Series A:
  Mathematical, Physical and Engineering Sciences 359~(1781) (2001) 739--757.

\bibitem{szaksz2023delay}
B.~Szaksz, G.~Stepan, Delay effects in the dynamics of human controlled towing
  of vehicles, Journal of Computational and Nonlinear Dynamics (2023) 1--12.

\bibitem{altintas2020chatter}
Y.~Altintas, G.~Stepan, E.~Budak, T.~Schmitz, Z.~M. Kilic, Chatter stability of
  machining operations, Journal of Manufacturing Science and Engineering
  142~(11) (2020).

\bibitem{stepan2014cylindrical}
G.~St{\'e}p{\'a}n, J.~Munoa, T.~Insperger, M.~Surico, D.~Bachrathy,
  Z.~Domb{\'o}v{\'a}ri, Cylindrical milling tools: Comparative real case study
  for process stability, CIRP Annals 63~(1) (2014) 385--388.

\bibitem{kiss2022control}
A.~K. Kiss, T.~G. Molnar, A.~D. Ames, G.~Orosz, Control barrier functionals:
  Safety-critical control for time delay systems, arXiv preprint
  arXiv:2206.08409 (2022).

\bibitem{hale2013introduction}
J.~K. Hale, S.~M.~V. Lunel, Introduction to functional differential equations,
  Vol.~99, Springer Science \& Business Media, 2013.

\bibitem{stepan1989retarded}
G.~St{\'e}p{\'a}n, Retarded dynamical systems: stability and characteristic
  functions, Longman Scientific \& Technical, 1989.

\bibitem{yoshida2021basins}
K.~Yoshida, K.~Konishi, N.~Hara, Basins and bifurcations of a delayed feedback
  control system and its experimental verification for a dc bus circuit,
  Nonlinear Dynamics 106 (2021) 2363--2376.

\bibitem{scholl2020norm}
T.~H. Scholl, V.~Hagenmeyer, L.~Gr{\"o}ll, On norm-based estimations for
  domains of attraction in nonlinear time-delay systems, Nonlinear Dynamics
  100~(3) (2020) 2027--2045.

\bibitem{michiels2007stability}
W.~Michiels, S.-I. Niculescu, Stability and stabilization of time-delay
  systems: an eigenvalue-based approach, SIAM, 2007.

\bibitem{leng2016basin}
S.~Leng, W.~Lin, J.~Kurths, Basin stability in delayed dynamics, Scientific
  reports 6~(1) (2016) 21449.

\bibitem{soliman1989integrity}
M.~Soliman, J.~Thompson, Integrity measures quantifying the erosion of smooth
  and fractal basins of attraction, Journal of Sound and Vibration 135~(3)
  (1989) 453--475.

\bibitem{thompson1989chaotic}
J.~Thompson, Chaotic phenomena triggering the escape from a potential well,
  Proceedings of the Royal Society of London. A. Mathematical and Physical
  Sciences 421 (1989) 195--225.

\bibitem{lenci2003optimal}
S.~Lenci, G.~Rega, Optimal control of homoclinic bifurcation: theoretical
  treatment and practical reduction of safe basin erosion in the helmholtz
  oscillator, Journal of Vibration and Control 9~(3-4) (2003) 281--315.

\bibitem{rega2005identifying}
G.~Rega, S.~Lenci, Identifying, evaluating, and controlling dynamical integrity
  measures in non-linear mechanical oscillators, Nonlinear Analysis: Theory,
  Methods \& Applications 63~(5-7) (2005) 902--914.

\bibitem{insperger2011semi}
T.~Insperger, G.~St{\'e}p{\'a}n, Semi-discretization for time-delay systems:
  stability and engineering applications, Vol. 178, Springer Science \&
  Business Media, 2011.

\bibitem{molnar2018bifurcation}
T.~G. Moln{\'a}r, Z.~Domb{\'o}v{\'a}ri, T.~Insperger, G.~St{\'e}p{\'a}n,
  Bifurcation analysis of nonlinear time-periodic time-delay systems via
  semidiscretization, International Journal for Numerical Methods in
  Engineering 115~(1) (2018) 57--74.

\bibitem{burton2005volterra}
T.~A. Burton, Volterra integral and differential equations, Vol. 202, Elsevier,
  2005.

\bibitem{burton2016existence}
T.~A. Burton, Existence and uniqueness results by progressive contractions for
  integro-differential equations, Nonlinear Dyn. Syst. Theory 16~(4) (2016)
  366--371.

\bibitem{vio2007bifurcation}
G.~A. Vio, G.~Dimitriadis, J.~E. Cooper, Bifurcation analysis and limit cycle
  oscillation amplitude prediction methods applied to the aeroelastic galloping
  problem, Journal of Fluids and Structures 23~(7) (2007) 983--1011.

\bibitem{cirillo2016spectral}
G.~I. Cirillo, A.~Mauroy, L.~Renson, G.~Kerschen, R.~Sepulchre, A spectral
  characterization of nonlinear normal modes, Journal of Sound and Vibration
  377 (2016) 284--301.

\bibitem{cenedese2022data}
M.~Cenedese, J.~Ax{\aa}s, B.~B{\"a}uerlein, K.~Avila, G.~Haller, Data-driven
  modeling and prediction of non-linearizable dynamics via spectral
  submanifolds, Nature communications 13~(1) (2022) 872.

\bibitem{insperger2002semi}
T.~Insperger, G.~St{\'e}p{\'a}n, Semi-discretization method for delayed
  systems, International Journal for numerical methods in engineering 55~(5)
  (2002) 503--518.

\bibitem{habib2017chatter}
G.~Habib, G.~Kerschen, G.~Stepan, Chatter mitigation using the nonlinear tuned
  vibration absorber, International Journal of Non-Linear Mechanics 91 (2017)
  103--112.

\bibitem{sims2007vibration}
N.~D. Sims, Vibration absorbers for chatter suppression: A new analytical
  tuning methodology, Journal of Sound and Vibration 301~(3-5) (2007) 592--607.

\bibitem{dombovari2008estimates}
Z.~Dombovari, R.~E. Wilson, G.~Stepan, Estimates of the bistable region in
  metal cutting, Proceedings of the Royal Society A: Mathematical, Physical and
  Engineering Sciences 464~(2100) (2008) 3255--3271.

\bibitem{engelborghs2002numerical}
K.~Engelborghs, T.~Luzyanina, D.~Roose, Numerical bifurcation analysis of delay
  differential equations using {DDE-BIFTOOL}, ACM Transactions on Mathematical
  Software (TOMS) 28~(1) (2002) 1--21.

\bibitem{belardinelli2016first}
P.~Belardinelli, S.~Lenci, A first parallel programming approach in basins of
  attraction computation, International Journal of Non-Linear Mechanics 80
  (2016) 76--81.

\bibitem{habib2019suppression}
G.~Habib, Suppression of time-delayed induced vibrations through the dynamic
  vibration absorber: application to the inverted pendulum, in: Topics in
  Nonlinear Mechanics and Physics: Selected Papers from CSNDD 2018, Springer,
  2019, pp. 125--140.

\bibitem{scarciotti2017nonlinear}
G.~Scarciotti, A.~Astolfi, et~al., Nonlinear model reduction by moment
  matching, Foundations and Trends{\textregistered} in Systems and Control
  4~(3-4) (2017) 224--409.

\bibitem{haller2016nonlinear}
G.~Haller, S.~Ponsioen, Nonlinear normal modes and spectral submanifolds:
  existence, uniqueness and use in model reduction, Nonlinear dynamics 86
  (2016) 1493--1534.

\bibitem{chaturantabut2010nonlinear}
S.~Chaturantabut, D.~C. Sorensen, Nonlinear model reduction via discrete
  empirical interpolation, SIAM Journal on Scientific Computing 32~(5) (2010)
  2737--2764.

\end{thebibliography}
	

\end{document}